\documentclass[matbio,10pt]{svjour}
\usepackage{inputenc}
\usepackage{amsmath}
\usepackage{amssymb,amsfonts}
\usepackage{graphicx,epstopdf,bm}
\usepackage{natbib}
\usepackage[colorlinks,bookmarks,urlcolor=blue,citecolor=blue,linkcolor=blue]{hyperref}
\usepackage[T1]{fontenc}
\usepackage{bbm} 
\usepackage{stmaryrd}

\newcommand{\pd}[2]{\frac{\partial #1}{\partial #2}}
\newcommand{\pdd}[2]{\frac{\partial^{2} #1}{\partial #2 ^{2}}}

\newcommand{\abs}[1]{\left\vert #1 \right\vert}
\newcommand{\ave}[1]{\left \langle #1 \right \rangle}
\newcommand{\explr}[1]{\exp\left[ #1 \right]}
\newcommand{\bigo}{O}
\newcommand{\diag}[1]{\Sigma_{#1}}
\newcommand{\chem}[2]{ {{\scriptscriptstyle#1\atop\longrightarrow}\atop{\longleftarrow\atop \scriptscriptstyle#2}} }
\newcommand{\bvec}[2]{\begin{bmatrix}#1\\#2 \end{bmatrix}}
\newcommand{\rmp}{\mathrm{p}}
\newcommand{\hgam}{\varphi}
\newcommand{\bp}{\mathbf{p}} 
\newcommand{\bv}{\mathbf{v}}
\newcommand{\bb}{\mathbf{b}} 
\newcommand{\br}{\mathbf{r}}
\newcommand{\bbe}{\mathbbm{e}}
\newcommand{\bk}{\mathbf{n}}
\newcommand{\bbW}{\mathbb{W}}
\newcommand{\bw}{\mathbf{w}}
\newcommand{\bu}{\mathbf{u}}
\newcommand{\anot}{\sigma}
\newcommand{\synthrate}{\tau}
\newcommand{\qd}{\alpha_{\rm e}}
\newcommand{\qss}{\alpha_{\rm i}}
\newcommand{\pvar}{p}
\newcommand{\tv}{y}
\newcommand{\bdrate}{\widehat{W}}
\newcommand{\hbdrate}{W}
\newcommand{\bpstat}{\hat{\mathbf{p}}} 
\newcommand{\isdt}{\mathbf{w}}
\newcommand{\esdt}{u}
\newcommand{\isd}{\hat{\isdt}} 
\newcommand{\esd}{\hat{u}} 
\newcommand{\bpss}{\bm{\rho}}  
\newcommand{\pss}[1]{\rho_{#1}}  
\newcommand{\bfone}{\mathbf{1}}  
\newcommand{\befn}[1]{\bm{\phi}_{#1}}  
\newcommand{\efn}[1]{\phi_{#1}}  
\newcommand{\baefn}[1]{\bm{\xi}_{#1}}  
\newcommand{\aefn}[1]{\xi_{#1}}  
\newcommand{\efnref}[1]{\efn{#1}^{(r)}}
\newcommand{\efnabs}[1]{\efn{#1}^{(a)}}
\newcommand{\aefnref}[1]{\aefn{#1}^{(r)}} 
\newcommand{\aefnabs}[1]{\aefn{#1}^{(a)}} 
\newcommand{\eigref}[1]{\lambda_{#1}^{(r)}} 
\newcommand{\eigabs}[1]{\lambda_{#1}^{(a)}} 
\newcommand{\mv}[1]{\ave{\bv(#1)}}  
\newcommand{\difu}{D}  

\numberwithin{equation}{section}

\begin{document}
\title{Metastable behavior in Markov processes with internal states}
\author{Jay Newby \and Jon Chapman}
\institute{Mathematical Institute, University of Oxford, 24-29 St Giles', Oxford, OX1 3LB, UK, \email{newby@maths.ox.ac.uk}}

\maketitle
\begin{abstract}
A perturbation framework is developed to analyze metastable behavior in stochastic processes with random internal and external states.
The process is assumed to be under weak noise conditions, and the case where the deterministic limit is bistable is considered.
A general analytical approximation is derived for the stationary probability density and the mean switching time between metastable states, which includes the pre exponential factor.
The results are illustrated with a model of gene expression that displays bistable switching.
In this model, the external state represents the number of protein molecules produced by a hypothetical gene.
Once produced, a protein is eventually degraded.
The internal state represents the activated or unactivated state of the gene; in the activated state the gene produces protein more rapidly than the unactivated state.
The gene is activated by a dimer of the protein it produces so that the activation rate depends on the current protein level.
This is a well studied model, and several model reductions and diffusion approximation methods are available to analyze its behavior.
However, it is unclear if these methods accurately approximate long-time metastable behavior (i.e., mean switching time between metastable states of the bistable system).
Diffusion approximations are generally known to fail in this regard.  
\end{abstract}
\section{Introduction}
A common feature found in many stochastic models of biological processes is a distinction between internal and external states \citep{vankampen79a}.
There are numerous examples of such Markov processes used as models for biological phenomena \citep{othmer88a,bicout97a,kepler01a,friedman05a,newby10b}.
Examples of an internal state include the number of open ion channels in the membrane of a neuron that affect its membrane voltage \citep{keener11a} and the on/off state of a gene that affects its protein production rate \citep{newby12a}.
The distinction between internal and external states should not be confused with the concept of intrinsic and extrinsic noise (see \citep{thattai01a} for an example related to gene expression).
A system with internal degrees of freedom is a classical idea in physics and applied mathematics, and the extension of this concept to Markov processes with internal states is well known in the literature \citep{hill85a,kramli83a,landman77a}.

Consider the following two stochastic processes: the discrete internal state, $S(t)$, and the external state, $X(t)$. 
We consider two possibilites: $X(t)\in \mathbb{Z}$ and $X(t)\in \mathbb{R}$ (i.e., a discrete jump process and a continuous process).
If $X(t)$ is independent of $S(t)$, there is one source of noise affecting $X(t)$, and we assume it is scaled by $1/\qd$, where $\qd\gg 1$, so that in the limit $\qd\to \infty$, $X(t)$ is a deterministic process.
Now consider the combined process where $S(t)$ and $X(t)$ are coupled.
In this case, there is a second source of noise affecting $X(t)$ through its dependence on $S(t)$.
We assume that there is a second large parameter, $\qss$, such that in the limit $\qss \to \infty$, the frequency of jumps in $S(t)$ becomes infinite.
In this limit, $X(t)$ depends only on the average value of $S(t)$, effectively eliminating the second noise source.
In the limit $\qss\to\infty$, $\qd \to \infty$, the combined process $(S(t), X(t)) \to \bar{x}(t)$, where $\bar{x}(t)\in \mathbb{R}$ is deterministic.

Under weak noise conditions, meaning close to the deterministic limit with $\qss \gg 1$ and $\qd \gg 1$, the dynamics of the deterministic system strongly influence the dynamics of the stochastic process.  In particular, we are interested in the case where the deterministic system has multiple stable solutions depending on the initial conditions.  On short timescales, a trajectory of the stochastic process fluctuates about the deterministic trajectory that has the same initial conditions.  However, metastable behavior in the stochastic process is not seen in the deterministic system because it depends on a small amount of noise present in the system to cause a transition from one of the stable deterministic solutions to the other.  Metastable transitions occur on a long timescale.

Metastable behavior is important because it represents fluctuation-induced phenomena not present in the deterministic system.
The standard example of metastable behavior is Brownian motion in a double well potential.
On short timescales the particle is most likely found near one of the two minima, and on long timescales the particle can transition over the energy barrier that separates each well.
Metastable transitions by nonlinear Markov processes with an internal and external state are more difficult to analyze than diffusion in a potential well, and exact analytical solutions are rarely possible.
Moreover, using Monte Carlo simulations to generate exact trajectories that display metastable behavior requires too much processor time to be practical.
It is therefore necessary to develop approximation methods.

One approximation method is to reduce the complexity of the model by eliminating a noise source.
Noise in the internal state is eliminated in the adiabatic limit, $\qss \to \infty$, where $S(t)$ is averaged out of $(S(t), X(t))$ to obtain a Markov process that approximates $X(t)$.
In other words, although $X(t)$ is not Markovian due to its dependence on $S(t)$, it may be approximately Markovian.
Eliminating noise in the external state with $\qd\to\infty$ results in a velocity jump process where the external state evolves deterministically in between random jumps in the internal state.
However, the timescale for a metastable transition is very sensitive to both the type of noise and the noise strength, and eliminating a noise source can lead to large errors.

Another way to reduce the complexity of the model is with a diffusion approximation obtained using a quasi-steady-state (QSS) reduction \citep{gardiner83a,thomas12a}.
This is very similar to the adiabatic limit, but uses a perturbation approach so that higher order terms can be included that account for noise in the internal state.
The QSS reduction also approximates $(S(t), X(t))$ with a single continuous Markov process for $X(t)$, but includes effects from both noise sources.
This approximation reduces the problem to diffusion in a double well potential.
The underlying assumption behind the QSS reduction is that $S(t)$ is well approximated by a random variable chosen from its steady-state distribution conditioned on a fixed value of $X(t)$.
While the QSS diffusion approximation is a useful tool in most circumstances, it is not accurate for characterizing metastability.

To describe metastable behavior, it is necessary to approximate both the effective potential and the timescales for metastable transitions.
For a 1D continuous Markov process on the state space $x\in \mathbb{R}$, the potential is straightforward to define, and if a stationary solution exists, it must have zero probability flux everywhere.
Given its usefulness at describing the qualitative features, we would like to know if we can define an effective potential in general.
For higher dimensional continuous Markov processes, the potential is no longer well defined when the curl of the drift velocity field is nonzero, and it is possible for the stationary density to exhibit a nonzero probability flux.
This is closely related to detailed balance conditions and thermodynamic equilibrium.
Developing a systematic formalism to describe nonequilibrium stationary behavior is particularly relevant in biology.  
It turns out that an effective potential can still be defined using perturbation theory \citep{schuss10a,ludwig75a,matkowsky83a,talkner87a,naeh90a,maier97a,hanggi84a,dykman94a} and large deviation theory \citep{freidlin98a,metzner09b,heymann08a}.
These tools can also be used to approximate the timescale associated with metastable transitions.
The methods presented here fit within the perturbation framework.

The theory of large deviations \citep{freidlin98a,shwartz95a,feng06a} is the mathematical foundation for the techniques used to study metastable transitions (rare events).  Here, we focus on perturbation-theory-based techniques \citep{schuss10a}, which we refer to as the quasi-stationary analysis (QSA).  
Large deviation theory provides rigorous results and error estimates, but does not provide a means of explicitly calculating the pre exponential factor (see Section \ref{sec:wkb}), which is part of the leading order transition time and stationary density approximations.
The QSA is formal but systematic and generally more practical for applications.
The QSA was developed to analyze the differential Chapman--Kolmogorov (CK) equation, which describes the process by its probability density function.  For a continuous Markov process, the QSA is well-developed for the Fokker--Planck equation \citep{schuss10a,ludwig75a,matkowsky83a,talkner87a,naeh90a,maier97a}.  The QSA has also be applied to the Master equation to analyze certain birth-death processes \citep{hanggi84a,dykman94a,hinch05a,doering05a,vellela07a,doering07a,escudero09a,bressloff10b,assaf11a}.  (The Fokker--Planck and Master equation are instances of the more general CK equation \citep{gardiner83a}.)
However, for weak noise problems where adiabatic elimination (i.e., stochastic averaging) of one noise source is necessary to reach the deterministic system, no one has developed these methods (as far as we are aware) to study metastable behavior without first applying a QSS-type diffusion approximation or adiabatic reduction, which has been shown to result in significant errors \citep{freidlin98a,newby10b}.
Recently, the QSA has been developed for the velocity-jump process (also sometimes called a piecewise deterministic process or hybrid process) \citep{keener11a,newby11b,newby12a}, which is the simplest example of a process with internal states.
In this paper, we further develop the QSA for the case where $S(t)$ and $X(t)$ are both intrinsically stochastic.

There are several advantages to the QSA.  First, if the process includes a discrete state, the QSA provides an approximation that accounts for all moments of the jump propagator (infinitesimal generator), whereas the diffusion approximations include only the first two moments (e.g., a diffusion approximation of a discrete jump Markov process by truncation of a Kramers--Moyal (KM) expansion).  Second, it provides a uniformly accurate approximation of the stationary probability density function.  Third, physically meaningful quantities, such as the effective potential and metastable transition rates, can be generalized to processes that do not assume detailed balance. 
Finally, the QSA can be applied to higher dimensional (by which we mean the deterministic limit $\bar{x}\in \mathbb{R}^{d}$, $d\geq 1$) nonequilibrium processes.

The main goal of the paper is to develop the QSA for a general class of Markov processes that have a discrete internal state, and we illustrate the analysis using a simple example problem.
The example problem is ideal because we can derive several approximations that serve as examples to which we apply the general QSA.
Since each is an approximation of a single model, we can compare the effects of metastability in different types of Markov processes.
In particular, we are interested in approximating two quantities: the timescales for metastable transitions and the effective potential.
The analysis of the example problem should inform our understanding about when reduction techniques, such as a diffusion approximation, fail to approximate these two quantities and why.
Previous work has shown that diffusion approximations lead to errors in both the $\qss \to \infty$ \citep{hanggi84a,walczak05a} and $\qd\to\infty$ limits \citep{newby12a}.  But what happens when both noise sources are present?  When is one noise source more significant than the other?  If the QSS reduction fails, why does it fail?   Is it due to large deviation errors like the system-size expansion, or is it because the QSS assumption is invalid?  Does it fail for the same reasons in each limit?

The paper is organized as follows. 
First, in Section \ref{sec:two-classes-markov}, we describe in detail two related versions of $(S(t), X(t))$: one where $X(t)$ is a discrete birth-death process and one where $X(t)$ is a continuos Markov process.
The quasi-stationary analysis is presented in Section \ref{sec:long-time-asymptotic}.
Then, in Section \ref{sec:model}, we introduce the example problem along with various approximations and model reductions.
After we apply the QSA to the example problem, results are presented in Section \ref{sec:results}.

\section{Two Markov processes with a discrete internal state}
\label{sec:two-classes-markov}
Consider the following coupled stochastic processes.
Let $S(t)$ be the discrete internal state on a finite state space having $M$ states, and let $X(t)$ be the external state.
The main conceptual difference between an internal and external state is that the dynamics of $S(t)$ is fast compared to $X(t)$, and in the deterministic limit, the effect of $S(t)$ on $X(t)$ is in some sense ``averaged out'' so that only $X(t)$ is observable.

To make the QSA as general as possible, we consider both a continuous and a discrete external state.
In the case of a discrete external state, we assume conditions under which a continuous approximation is valid.
The joint probability density function (probability mass function if $X(t)$ is discrete) can be written as
\begin{align*}
   \rmp(s, x, t)\Delta x &= \text{Pr}[S(t)=s, X(t)\in (x, x+\Delta x)] \\
   &= \text{Pr}[ S(t)=s |X(t)\in (x, x+\Delta x)]  \text{Pr}[X(t)\in (x, x+\Delta x)].
\end{align*}
Define the conditional internal state distribution to be
\begin{equation}
  \label{eq:59}
  \mathrm{w}(s, t | x) \equiv \text{Prob}[ s=S(t) | X(t) = x],
\end{equation}
and the marginal external state density function to be
\begin{equation}
  \label{eq:60}
  \esdt(x, t) \equiv \text{Prob}[X(t)\in (x, x+\Delta x)]/\Delta x.
\end{equation}
It is convenient to use vector notation for the probability density with
\begin{equation}
  \label{eq:48}
  \mathbf{p}(x,t) \equiv (\mathrm{p}(0, x,t), \mathrm{p}(1, x,t), \cdots, \mathrm{p}(M-1, x,t))^{T}.
\end{equation}
In general, we sometimes represent a given function $f(s)$ as the vector $\mathbf{f}\in \mathbb{R}^{M}$ where the $s$th component of ${\bf f}$ is $f(s)$.  Diagonal matrices are written as $\diag{\bf f}$, where the diagonal entries are given by the elements of the vector $\mathbf{f}$,  and occasionally we may use the notation $\diag{f(s)}$, defined as $\diag{f(s)} = \diag{\bf f}$.

We write the stationary versions of \eqref{eq:59}-\eqref{eq:48} as
\begin{equation}
  \label{eq:41}
  \lim_{t\to\infty}\mathbf{p}(x, t) = \bpstat(x) =  \isd(x) \esd(x),
\end{equation}
where $\isd(x) \in \mathbb{R}^{M}$ and 
\begin{equation}
  \label{eq:42}
  \sum_{s}\hat{w}(s | x) = 1.
\end{equation}

For a fixed external state $X(t) = x$, the process $S(t)$ is described by a Master equation
\begin{equation}
  \label{eq:18}
  \pd{\mathbf{w}}{t} = \qss A(x)\mathbf{w}(t | x),
\end{equation}
where $A$ is a transition rate matrix and $\qss\gg 1$ is a large parameter.
The matrix $A$ is a member of a family of matrices called $\bbW$-matrices, which have the following properties.
First, the columns sum to zero, which means that the matrix is singular and the vector, $\bfone \equiv (1, \cdots, 1)^{T}$, is the left eigenvector corresponding to a zero eigenvalue.
For a transition rate matrix to be a $\bbW$-matrix, it must have negative diagonal elements, nonnegative off-diagonal elements, and it must be irreducible.
One can show, using the {\em Perron-Frobenius theorem}, that the nullspace of a $\bbW$-matrix is one dimensional and that the right nullvector has strictly positive elements.
Hence, there exists a unique $\bpss >0$ such that 
\begin{equation}
  \label{eq:19}
  A\bpss=0,\quad \sum_{s}\rho(s | x) = 1.
\end{equation}
For a fixed external state, $\rho(s | x)$ is the steady state distribution of the internal state, and we refer to it as the quasi-steady-state distribution.

We call the coupled process, $(S(t), X(t))$, with $X(t)$ continuous the semi-continuous process.
In this case the external state is given by the Ito stochastic differential equation,
\begin{equation}
  \label{eq:66}
  dX(t) = -v(S(t), X(t))dt + \sqrt{\frac{b(S(t), X(t))}{\qd}}dW(t),
\end{equation}
where $dW(t)$ is a Wiener process, $v(s, x)$ is the drift, $b(s, x)$ is the scaled diffusivity, and $\qd\gg 1$ is a large parameter.
The coupled process, $(S(t), X(t))$, is described by the CK equation
\begin{equation}
  \label{eq:52}
  \pd{}{t}\mathbf{p}(x, t) = \qss A(x)\mathbf{p} + \diag{\mathbbm{g}(s)}\mathbf{p} ,
\end{equation}
where the operator, $\mathbbm{g}$, is defined by
\begin{equation}
  \label{eq:68}
  \mathbbm{g}(s)\mathrm{p}(s, x, t) \equiv - \pd{}{x}(v(s, x)\mathrm{p}) + \frac{1}{2\qd}\pdd{}{x}\left(b(s, x) \mathrm{p}\right).
\end{equation}

The coupled process with $X(t)$ discrete is referred to as the discrete process.
The external state is defined in terms of the birth/death process, $N(t)\in \mathbb{Z}_{+}$, satisfying
\begin{equation}
  \label{eq:20}
  N(t) = N(0) + Y_{+}(\int_{0}^{t}\bdrate_{+}(N(\tau)|S(\tau))d\tau) 
  - Y_{-}(\int_{0}^{t}\bdrate_{-}(N(\tau)|S(\tau))d\tau),
\end{equation}
where $\bdrate_{+}$ and $\bdrate_{-}$ are the birth and death rates, respectively, and $Y_{\pm}(t)$ are unit Poisson processes.
We assume that the rates can be written as $\bdrate_{\pm}(n|s) = \qd\hbdrate_{\pm}(n/\qd|s)$ where $\qd \gg 1$ is a large parameter.  Let $X(t) = N(t)/\qd$.  
Then, \eqref{eq:20} can be written as
\begin{multline}
  \label{eq:45}
  X(t) = X(0) + \frac{1}{\qd}Y_{+}(\qd\int_{0}^{t}\hbdrate_{+}(X(\tau)|S(\tau))d\tau) \\
- \frac{1}{\qd}Y_{-}(\qd\int_{0}^{t}\hbdrate_{-}(X(\tau)|S(\tau))d\tau).
\end{multline}
The CK equation describing the discrete process is 
\begin{equation}
  \label{eq:53}
  \pd{}{t}\mathbf{p}(x, t) =  \qd A(x)\mathbf{p} + \diag{\mathbbm{d}(s)} \mathbf{p},
\end{equation}
where the operator, $\mathbbm{d}$, is defined by
\begin{equation}
  \label{eq:67}
   \mathbbm{d}(s)\mathrm{p}(s, x, t) \equiv \qd\left[(\bbe^{\partial x} - 1)\hbdrate_{-}(x|s)\mathrm{p} + (\bbe^{-\partial x} - 1)\hbdrate_{+}(x|s)\mathrm{p}\right].
\end{equation}
The jump operator,
\begin{equation}
  \label{eq:55}
  \bbe^{\pm\partial x}f(x) \equiv f(x\pm \frac{1}{\qd}) = \sum_{n=0}^{\infty}\frac{(\pm 1)^{n}}{\qd^{n}n!}f^{(n)}(x),
\end{equation}
can be written in terms of a Taylor series expansion, which formally yields the Kramers--Moyal expansion of \eqref{eq:53}.

Note that if $v(s, x)$ and $b(s, x)$ are chosen appropriately, the semi-continuous process \eqref{eq:66} is a diffusion approximation of the birth death process \eqref{eq:45}.
Without loss of generality, we assume that the two external state processes are related by
\begin{equation}
  \label{eq:117}
  v(s, x) = \hbdrate_{+}(x | s) - \hbdrate_{-}(x | s), \quad b(s, x) = \frac{1}{2}\left(\hbdrate_{+}(x | s) + \hbdrate_{-}(x | s)\right).
\end{equation}

Metastable behavior requires the stochastic processes to be under weak noise conditions.
Generally speaking, the QSA is an asymptotic analysis where a small variable, call it $\epsilon$, controls the global noise strength.
That is, in the limit $\epsilon \to 0$, the stochastic process converges to a deterministic system.
As discussed in the Introduction, the purpose of the two large parameters $\qss$ and $\qd$ is to place the stochastic process in weak noise conditions.
In the limit $\qss\to \infty$, noise from the internal state is eliminated.
Likewise, in the limit $\qd\to\infty$, noise from the external state is eliminated.
Hence, both limits must be taken to reach a deterministic system.
In order to carry out a systematic asymptotic analysis with a single small parameter, we define $\epsilon = 1/\qss = 1/(\hgam\qd)$.
Hence, the limit $\epsilon\to 0$ is equivalent to taking the limit $\qss \to \infty$, $\qd\to\infty$ with the ratio $ \hgam = \qss/\qd$ fixed.

We assume for either process that the deterministic limit,
\begin{equation}
  \label{eq:114}
  \dot{x} = \bar{v}(x) \equiv \sum_{s}\rho(s | x)v(s, x),
\end{equation}
is bistable.  That is, there are three fixed points satisfying $\bar{v}(x) = 0$, label them $x_{-} < x_{*} < x_{+}$, with $\bar{v}(x)>0$ for $x<x_{-}$, $\bar{v}(x)<0$ for $x_{-} < x < x_{*}$, $\bar{v}(x) > 0$ for $x_{*} < x < x_{+}$, and $\bar{v}(x) < 0$ for $x > x_{+}$.  Then, $x_{\pm}$ are stable fixed points and $x_{*}$ is unstable.

To ensure a well-defined process, we assume for some interval $(x_{a}, x_{b})$, with $x_{a} < x_{-} < x_{+} < x_{b}$, that $v$, $b$, and $\hbdrate_{\pm}$ are smooth functions of $x$.  Assume further that $\hbdrate_{\pm}(x | s) > 0$,  for all $0 \leq s \leq M-1$ and $x \in (x_{a}, x_{b})$.\footnote{The last constraint can be relaxed somewhat provided the process converges to a unique stationary density.}

\section{Quasi-stationary analysis}
\label{sec:long-time-asymptotic}
\newcommand{\intst}{s}
We now present a systematic perturbation method to analyze metastable, or long-time, behavior of the discrete and semi-continuous processes.
Suppose we have a CK equation of the form
\begin{equation}
  \label{eq:32}
  \pd{}{t}\mathrm{p}(s, x, t)  = -\mathcal{L}_{\epsilon}\mathrm{p},
\end{equation}
where $\mathcal{L}_{\epsilon}$ is a compact linear operator acting on functions of $(s, x)$.  Note that one can easily generalize this theory to the case $x\in \Omega\subset\mathbb{R}^{N}$ (see \citep{newby12a}).   For illustration, take $\mathcal{L}_{\epsilon}$ to have the form
\begin{multline}
  \label{eq:94}
  \mathcal{L}_{\epsilon}\mathrm{p} \equiv \pd{}{x}(v(s, x)\mathrm{p}) - \epsilon\pdd{}{x}(b(s, x)\mathrm{p}) \\
- \frac{1}{\epsilon}\sum_{\intst'}\left(A(s, s' | x)\mathrm{p}(s', x, t) - A(s', s | x)\mathrm{p}(s, x, t)\right),
\end{multline}
where $\epsilon\ll 1$ is a small parameter.
(Note that we have absorbed $\hgam$ into the definition of $b(s, x)$.)

Assume that $\mathcal{L}_{\epsilon}$ has a complete set of eigenfunctions, $\{\efn{j}(s, x)\}$ and adjoint eigenfunctions $\{\aefn{j}(s, x)\}$.  If the initial condition is $\mathrm{p}(s, x, 0) = \delta(x-x_{0})\delta_{s,s_{0}}$, the solution can be written
\begin{equation}
  \label{eq:33}
  \mathrm{p}(s, x, t) = \sum_{j=0}^{\infty}\aefn{j}(s_{0}, x_{0})\efn{j}(s, x)e^{-\lambda_{j}t},
\end{equation}
where we assume that all of the eigenvalues, $\lambda_{j}$, are nonnegative.  Since we are interested in metastable behavior, assume that in the limit $\epsilon\to 0$, $X(t)$ converges to a bistable deterministic process.  Label the two stable fixed points $x_{\pm}$ and the unstable fixed point $x_{*}$ and assume $x_{-}<x_{*}<x_{+}$.

The random process will look very different if the external state starts at $x_{0}<x_{*}$ or $x_{0}>x_{*}$.  For the sake of illustration assume that $x_{0}=x_{-}$.  On intermediate time scales, the solution will converge to a stationary density around $x_{-}$ that, figuratively speaking, does not see the other stable fixed point---or said another way, the solution does not see beyond $x_{*}$.  Slowly, over a long timescale, the solution converges to the full stationary density as probability slowly leaks out past $x_{*}$ toward $x_{+}$.  The timescale for this long-time convergence is exponentially large (i.e., $O(e^{1/\epsilon}$)).  Since a stationary solution exists, the smallest eigenvalue $\lambda_{0}$, called the principal eigenvalue, is $\lambda_{0}=0$, and the stationary density is the eigenfunction $\efn{0}(s, x)$; that is, we normalize the principal eigenfunction so that $\ave{\efn{0}, 1}=1$ with respect to the inner product defined by
\begin{equation}
  \label{eq:203}
  \ave{f(s, x), g(s, x)} \equiv \int_{x\in \Omega}\sum_{\intst}f(s, x)g(s, x)dx.
\end{equation}

The separation of time scales in the problem can be exploited to approximate the solution.  To understand how this works consider the process where a boundary condition is placed at $x_{*}$ so that the process truly does not see beyond the unstable fixed point.  We want to consider two different boundary conditions: reflecting and absorbing.  
To distinguish between each case, we write the principal eigenvalue and eigenfunction (dropping the subscript) as $\eigabs{},\;\efnabs{}$ and $\eigref{},\;\efnref{}$ for absorbing and reflecting boundary conditions, respectively.  If we place a reflecting boundary at $x_{*}$ the principal eigenvalue $\eigref{}=0$, but the eigenfunction $\efnref{}$ is now restricted to $x\in \Omega_{-} = (-\infty, x_{*})$ (or $x\in \Omega_{+} = ( x_{*}, \infty)$ if we instead assume that $x_{0}>x_{*}$).  We call this the quasi-stationary density.  

Now suppose that an absorbing boundary is imposed at $x_{*}$.  In this case, no stationary density exists, and the principal eigenvalue is perturbed by an exponentially small amount, that is, $\eigabs{}=\bigo(e^{-C/\epsilon})$, for some $C>0$.  The eigenfunction $\efnabs{}$ is also perturbed, but away from the boundary, $\efnabs{} \sim \efnref{}$, which turns out to be straight forward to compute using a Wentzel--Kramers--Brillouin (WKB) approximation method.  Thus, if we can calculate the eigenvalue and eigenfunction, we have an accurate approximation to the absorbing boundary problem with
\begin{equation}
  \label{eq:34}
  \mathrm{p}(s, x, t) \sim \efnabs{}(s, x)e^{-\eigabs{} t}, \quad t\eigabs{1}\gg 1,
\end{equation}
or, since $\efnabs{}\sim \efnref{}$,
\begin{equation}
\label{eq:61}
  \mathrm{p}(s, x, t) \sim \efnref{}(s, x)e^{-\eigabs{} t}, \quad t\eigabs{1}\gg 1,\quad \epsilon\ll 1.
\end{equation}
We discuss how to approximate $\eigabs{}$ later in this section.  

This approximation can be repeated for the initial condition $x_{0}>x_{+}$, and a different principal eigenvalue and quasi-stationary density are obtained, call the eigenvalues $\eigabs{\pm}$ and quasi-stationary densities $\efnref{\pm}(s, x)$.  The full system, without any boundary condition imposed at $x_{*}$, can then be approximated by
\begin{equation}
  \label{eq:35}
  \mathrm{p}(s, x, t) \sim 
  \begin{cases}
    q_{-}(t)\efnref{-}(s, x), & x<x_{*} \\
    q_{+}(t) \efnref{+}(s, x),& x>x_{*}
  \end{cases},
\end{equation}
where $q_{\pm}(t)$ satisfy the system of ordinary differential equations
\begin{align}
  \label{eq:36}
  \frac{d q_{-}}{dt} &= -\eigabs{-}q_{-} + \eigabs{+}q_{+} \\
  \frac{d q_{+}}{dt} &= \eigabs{-}q_{-} - \eigabs{+}q_{+},
\end{align}
with $q_{-}(0)=1$ and $q_{+}(0)=0$ if $x_{0}<x_{*}$, or $q_{-}(0)=0$ and $q_{+}(0)=1$ if $x_{0}>x_{*}$.

A closely related problem is the mean escape time from a potential well.
Define the escape time, $\tau_{\pm}$, as the first time the process reaches $x_{*}$ having started at $x_{0} = x_{\pm}$.
Define the mean escape time as $T_{\pm} = \ave{\tau_{\pm}}$.
It follows from \eqref{eq:61} that $\tau_{\pm}$ can be approximated by an exponential random variable with mean $T_{\pm} \sim 1/\eigabs{\pm}$.

To obtain an approximation of the principal eigenvalues for each well, $\eigabs{\pm}$, we use a spectral projection method that makes use of the adjoint operator $\mathcal{L}_{\epsilon}^{*}$. 
(For simplicity, we drop the $\pm$ notation as the following analysis applies for either potential well.)
The spectral projection method was first developed for scalar-value PDE eigenvalue problems \citep{lee95a,hinch05a} and later generalized to a vector-valued PDE eigenvalue problem \citep{newby11b,keener11a,newby12a}.  The present treatment further generalizes the method. Consider the adjoint eigenfunctions $\{\aefn{j}\}$, $j=0,1,\cdots$, satisfying $\mathcal{L}_{\epsilon}^{*}\aefn{j} = \lambda_{j} \aefn{j}$,
and take $\ave{\efn{i},\aefn{j}} = \delta_{i,j}$ so that the two sets of eigenfunctions are biorthogonal.  We use the same notation to distinguish between the two boundary conditions for the adjoint eigenfunction.  If the boundary is reflecting, the first adjoint eigenfunction is $\aefnref{} = 1$, and if the boundary is absorbing then $\aefnabs{} \sim \aefnref{}$ away from the boundary, but develops a boundary layer at $x_{*}$. 
Using integration by parts we have
\begin{equation}
  \label{eq:39}
\ave{\efnref{},\eigabs{}\aefnabs{}} =  \ave{\efnref{},\mathcal{L}_{\epsilon}^{*}\aefnabs{}} = \ave{\mathcal{L}_{\epsilon} \efnref{},\aefnabs{}}+J(\efnref{},\aefnabs{}),
\end{equation}
where the boundary contribution,
\begin{equation}
  \label{eq:204}
  J(\efnref{}, \aefnabs{}) = \epsilon\sum_{\intst}b(s, x_{*})\efnref{}(s, x_{*})\frac{d}{dx}\aefnabs{}(s, x_{*}),
\end{equation}
is nonzero because $\efnref{}$ does not satisfy the absorbing boundary condition.  Then, since $\mathcal{L}_{\epsilon}\efnref{}=0$, the principal eigenvalue is
\begin{equation}  
   \label{eq:40}
  \eigabs{} = \frac{J(\efnref{},\aefnabs{})}{\ave{\efnref{},\aefnabs{}}}.
\end{equation}
The above identity can be used to approximate the principal eigenvalue as follows.  Since away from the boundary $x_{*}$, $\aefnabs{} \sim \aefnref{} = 1$, we can make this substitution for the term in the denominator of \eqref{eq:40} so that 
\begin{equation}
  \label{eq:38}
  \eigabs{} \sim \frac{J(\efnref{},\aefnabs{})}{\ave{\efnref{},1}},
\end{equation}
with exponentially small error.  Notice that the denominator is then well approximated by the normalization factor for the eigenfunction $\efnref{}$.  We cannot make the same substitution for the term in the numerator, since \eqref{eq:40} becomes a formula for $\eigref{}=0$ instead of $\eigabs{}$.  Of course, in some sense zero is actually a very good approximation because the error is $O(e^{-C/\epsilon})$, but to capture the metastable behavior we need to capture the small exponential.  Note that we could have just as well used $\efnabs{}$ and $\aefnref{}$ in \eqref{eq:39} instead of $\efnref{}$ and $\aefnabs{}$.  We choose the later because it simplifies the boundary layer analysis.

The recipe for approximating the solution requires approximations of the first eigenfunction and the first adjoint eigenfunction, where the latter satisfies the appropriate adjoint absorbing boundary condition.  In the remainder of this section, we calculate asymptotic approximations for the two eigenfunctions and then use the results to obtain an asymptotic approximation of $\eigabs{}$.  The main results are stated in Theorems \ref{pr:ef}-\ref{pr:lambda}.

\subsection{WKB approximation of the eigenfunction $\efnref{}(s, x)$}
\label{sec:wkb}

To simplify notation, we refer to $\efnref{}(s, x)$ as $\efn{}(s, x)$, and consistent with the vector notation introduced in Section \ref{sec:two-classes-markov}, we define the vector $\befn{}(x)$ as having elements given by $\efn{}(s, x)$, $s=0,1,\cdots,M-1$.
From the CK equations \eqref{eq:52} and \eqref{eq:53} it follows that the eigenfunction (up to terms exponentially small in $\epsilon$) satisfies
\begin{equation}
  \label{eq:228}
  \left[ A(x) + \frac{1}{\qss}\diag{\mathbbm{w}(s)}\right]\befn{}(x) = 0,
\end{equation}
where $\mathbbm{w} = \mathbbm{g}$ ($\mathbbm{w} = \mathbbm{d}$) for the semi-continuous (discrete) process with $\mathbbm{g}$ and $\mathbbm{d}$ defined by \eqref{eq:67} and \eqref{eq:68}, respectively.
We assume that the eigenfunction has the following WKB form
\begin{equation}
  \label{eq:201}
    \befn{}(x) \sim (\br_{0}(x)+\epsilon\br_{1}(x) + \cdots)\explr{-\frac{1}{\epsilon}\Phi(x)},
\end{equation}
where $\Phi$ is a scalar functions and $\br_{0,1}\in \mathbb{R}^{M}$  (with $\br_{0}$ positive).  Substituting \eqref{eq:201} into \eqref{eq:228} and collecting leading order terms yields
\begin{equation}
  \label{eq:227}
\bigo(1):\quad \left[A(x) + \diag{\mathbf{h}(x, \Phi'(x))}\right]\br_{0}(x) = 0,
\end{equation}
where,
\begin{equation}
  \label{eq:232}
  \mathrm{h}^{\mathrm{disc}}(s, x, \pvar) = \frac{1}{\hgam}\left[\hbdrate_{+}(x|s)(e^{-\hgam \pvar} - 1) + \hbdrate_{-}(x|s)(e^{\hgam \pvar} - 1)\right],
\end{equation}
for the discrete process and 
\begin{equation}
  \label{eq:231}
  \mathrm{h}^{\mathrm{sc}}(s, x, \pvar) = \pvar v(s, x) + \pvar^{2}b(s, x)
\end{equation}
for the semi-continuous process.  
For notational convenience, we have set $p = \Phi'$.
We rewrite the remaining term in \eqref{eq:227} as $\mathbf{r}_{0}(x) = k(x)\isd(x)$, where the approximation for the conditional internal state distribution \eqref{eq:59} is determined by calculating the nullspace of $A+\diag{\mathbf{h}}$.
Note that at fixed points, $x_{c}= x_{\pm},x_{*}$, we have that $\isd(x_{c}) = \bpss(x_{c})$, where $\bpss(x)$ is the quasi-steady-state distribution satisfying $A\bpss = 0$.
The scalar function $k(x)$ is a normalization factor, often referred to as the pre exponential factor in the literature, and is determined at higher order.
An equation for $\Phi'$ is given by
\begin{equation}
  \label{eq:58}
  \mathcal{H}(x, p) \equiv \det(A(x) + \diag{\mathbf{h}(x, p)}) = 0,
\end{equation}
where the function $\mathcal{H}(x, p)$ is called the Hamiltonian.
Since we must have $\mathbf{r}_{0} > 0$, a suitable solution to \eqref{eq:58} must result in a positive nullspace of $A(x) +  \diag{\mathbf{h}(x, p)}$.

To calculate $k(x)$, substitute \eqref{eq:201} into \eqref{eq:228} and collect $\bigo(\epsilon)$ terms to get
\begin{multline}
\nonumber
\bigo(\epsilon):  \quad \left[A + \diag{\mathbf{h}(x, \Phi'(x))}\right]\mathbf{r}_{1} = \frac{dk}{dx}\diag{\mathbf{h}_{\pvar}}\isd + k\diag{\mathbf{h}_{\pvar}}\frac{d\isd}{dx}\\
+ k \left(\diag{\mathbf{h}_{\pvar x}}+ \frac{1}{2}\Phi''(x)\diag{\mathbf{h}_{\pvar\pvar}}\right)\isd. 
\end{multline}
While the expansion is straightforward for the semi-continuous process, it is somewhat more complicated for the discrete process.  We leave the details to Appendix \ref{sec:wkbkm-expansion}.
We can use a solvability condition to derive an equation for $k(x)$ as follows.
Define the left nullvector, $\bm{l}$, with $ \bm{l}^{T}[A + \diag{\mathbf{h}(x, \Phi'(x))}] = 0$.
It follows from the Fredholm Alternative Theorem that $\mathbf{r}_{1}$ exists if and only if $k(x)$ satisfies
\begin{equation}
  \label{eq:245}
  \frac{dk}{dx} + \Psi'(x) k = 0,
\end{equation}
where
\begin{equation}
  \label{eq:242}
  \Psi'(x) = \frac{\bm{l}^{T}(x)\mathbf{H}_{\pvar x}(x, \Phi'(x)) + \frac{1}{2}\Phi''(x)\bm{l}^{T}(x)\mathbf{H}_{\pvar\pvar}(x, \Phi'(x))}{\bm{l}^{T}(x)\mathbf{H}_{\pvar}(x, \Phi'(x))},\quad x\neq x_{\pm},x_{*},
\end{equation}
with
\begin{equation}
  \label{eq:243}
  \mathbf{H}(x, \pvar) \equiv  [A + \diag{\mathbf{h}(x, \pvar)}]  \isd(x).
\end{equation}
 We can express $\Phi''(x)$ in terms of partial derivatives of the Hamiltonian with
\begin{equation}
  \label{eq:248}
  \Phi''(x) = -\frac{\mathcal{H}_{x}(x, \Phi'(x))}{\mathcal{H}_{\pvar}(x, \Phi'(x))}, \quad x\neq x_{\pm}, x_{*}.
\end{equation}
(For more about evaluating the limit $x\to x_{c}$, $x_{c} = x_{\pm},x_{*}$, of $\Phi''(x)$ and $\Psi'(x)$, see Appendix \ref{sec:fplim}.)
Hence,
\begin{equation}
  \label{eq:62}
  k(x) = \explr{-\Psi(x)}.
\end{equation}
\begin{theorem}
  \label{pr:ef}
Given a solution $p=\Phi'(x)$ to \eqref{eq:58} and its integral $\Phi(x)$, an asymptotic approximation of the solution to \eqref{eq:228} is given by
\begin{equation}
  \label{eq:63}
  \befn{}(x) \sim \isd(x)\explr{-\frac{1}{\epsilon}\Phi(x) - \Psi(x)},
\end{equation} 
where $\Psi(x)$ is given by integration of \eqref{eq:242}, $\mathbf{r}_{0} = k(x)\isd(x)$ satisfies \eqref{eq:227}, and $\isd(x)$ satisfies \eqref{eq:42}.
\end{theorem}

\subsection{Singular perturbation approximation of the adjoint eigenfunction $\aefnabs{}(s, x)$}
\label{sec:sing-pert-appr}
\newcommand{\meig}{\gamma}
The WKB method used in the previous section provides only an approximation of the stationary density, not the timescale for metastable transitions (i.e., the principal eigenvalue $\eigabs{}$).
To get information about transition times we must calculate an approximation of the adjoint eigenfunction.
As in Section \ref{sec:wkb}, we simplify notation with $\aefn{}(s, x) = \aefnabs{}(s, x)$ and define the vector $\baefn{}(x)$, having elements $\aefn{}(s, x)$, $s=0, 1, \cdots M-1$.
The analysis for the semi-continuous and discrete processes are sufficiently different that we present each separately.
\subsubsection{Semi-continuous process}
\label{sec:semi-cont-proc}
Up to terms exponentially small in $\epsilon$, the first adjoint eigenfunction satisfies
\begin{equation}
  \label{eq:86}
  \left[[A(x)]^{T} + \epsilon\diag{\bv(x)}\frac{d}{dx} + \epsilon^{2}\diag{\bb(x)}\frac{d^{2}}{dx^{2}}\right]\baefn{}(x) = 0,
\end{equation}
along with the absorbing boundary condition,
\begin{equation}
  \label{eq:87}
  \baefn{}(x_{*}) = 0.
\end{equation}
The outer solution, which does not satisfy the boundary condition, is exactly $\baefn{\mathrm{out}}\equiv \bfone$.

To obtain an approximate solution that also satisfies boundary conditions, we must rescale $x = x_{*}+\epsilon^{\theta}z$, for some $\theta>0$.  A reasonable first try is to take $\theta=1$ so that $x = x_{*}+\epsilon z$.  Equation \eqref{eq:86} becomes
\begin{equation}
  \label{eq:88}
   \left[[A(x_{*})]^{T} + \diag{\bv(x_{*})}\frac{d}{dz} +  \diag{\bb(x_{*})}\frac{d^{2}}{dz^{2}}\right]\baefn{\mathrm{bl}}(z) = 0,
\end{equation}
where $\baefn{\mathrm{bl}}(z) \equiv \baefn{}(x_{*} + \epsilon z)$.
The solution is a linear combination of the subsolutions $c_{j}\Upsilon_{j}e^{-\meig_{j}z}$, $j=0,\cdots,2M-1$ where
\begin{equation}
  \label{eq:89}
  \left( A(x_{*})^{T} - \meig_{j}\diag{\bv(x_{*})} + \meig_{j}^{2}\diag{\bb(x_{*})}\right)\Upsilon_{j} = 0,
\end{equation}
and $c_{j}$, $j=0,\cdots,2M-1$ are unknown constants.  We can specify the first solution as $\Upsilon_{0}=\bfone$ and $\meig_{0}=0$.  A valid solution should be bounded in the limit $z\to\infty$, which means that $c_{j}=0$  if  $\meig_{j}<0$; although we do not know a priori how many of the eigenvalues are negative.  Note that the boundary condition \eqref{eq:87} provides a system of $M$ linear equations for the $2M$ unknowns, $c_{j}$, which means that constraints to eliminate the remaining $M$ unknowns are required to close the system.  One such constraint eliminates an unknown (i.e., $c_{0}$) by matching to the outer solution, leaving $M-1$ more constraints we must find.  
We assume that there are $M-1$ negative eigenvalues.  For simplicity, we order the eigenvalues so that $\meig_{j}<0$ for $j = M+1, \cdots, 2M-1$.  

For the moment, consider the matrices in \eqref{eq:89} as depending on $x$ so that $\Upsilon_{j}$ and $\meig_{j}$ are also functions of $x$.  It is simple to show that $\Upsilon_{0}(x)=\bfone$ and $\meig_{0}(x)=0$ even if the matrices are evaluated away from $x_{*}$.  However, one of the solutions, label it $j=1$, is $\meig_{1} = \pvar(x) = \Phi'(x)$ from \eqref{eq:227}.  Moreover, $\Upsilon_{1}(x)\to \bfone$ and $\meig_{1}(x)\to 0$ as $x\to x_{*}$.  In fact, we know that $\meig_{1}(x)$ vanishes at all of the deterministic fixed points because $\Phi'(x_{c})=0$, for $x_{c}=x_{\pm},x_{*}$, and is nonzero otherwise.  It follows that the zero eigenvalue has a degenerate eigenspace, and the solution must include a secular term involving the generalized eigenvector satisfying
\begin{equation}
  \label{eq:90}
  A(x_{*})^{T}\bm{\zeta} =  \diag{\bv(x_{*})}\bfone=\bv(x_{*}).
\end{equation}
 One can show\citep{newby11b} that the deterministic fixed points are the only points where the eigenspace associated with the zero eigenvalue is degenerate.  The solution to \eqref{eq:88} is thus
\begin{equation}
  \label{eq:91}
    \baefn{\mathrm{bl}}(z) = c_{0}\bfone + c_{1}(\bm{\zeta} - z\bfone) + \sum_{j=2}^{M}c_{j}\Upsilon_{j}e^{-\meig_{j}z}.
\end{equation}
However, because of the secular term, the solution is unbounded in the limit $z\to\infty$, and as a result, it cannot be matched to the outer solution.  Therefore, there is a transition layer that sits between the boundary layer and the outer region.

To find the scaling for this transition layer, we change variables to $x = x_{*} + \epsilon^{\theta}\tv$, for $0<\theta<1$, and define $\baefn{\theta}(\tv) \equiv \baefn{}(x_{*} + \epsilon^{\theta} \tv)$.  Introduce the asymptotic expansion
\begin{equation}
  \label{eq:92}
   \baefn{\theta}(\tv) \sim \baefn{\theta}^{(0)}(\tv) + \epsilon^{\kappa} \baefn{\theta}^{(1)}(\tv) + \epsilon^{2\kappa} \baefn{\theta}^{(2)}(\tv),
\end{equation}
where $\kappa>0$.
Equation \eqref{eq:86} becomes
\begin{multline}
  \label{eq:93}
   \bigg[\left(A(x_{*}) + \epsilon^{\theta}\tv A'(x_{*}) + \cdots\right)^{T} \\
   + \epsilon^{1-\theta}\left(\diag{\bv(x_{*})} + \epsilon^{\theta}\tv\diag{\bv'(x_{*})} +\cdots\right)\frac{d}{d\tv} \\
\shoveright{ \left. + \epsilon^{2(1-\theta)}\left(\diag{\bb(x_{*})} +\cdots\right)\frac{d^{2}}{d\tv^{2}} \right]}\\
\times\left(\baefn{\theta}^{(0)}(\tv) + \epsilon^{\kappa} \baefn{\theta}^{(1)}(\tv) + \epsilon^{2\kappa} \baefn{\theta}^{(2)}(\tv)  \right) = 0.
\end{multline}
Setting $\epsilon=0$ in \eqref{eq:93} yields
\begin{equation}
  \label{eq:95}
  \bigo(1):\quad A(x_{*})^{T}\baefn{\theta}^{(0)}(\tv) = 0,
\end{equation}
which implies that
\begin{equation}
  \label{eq:96}
  \baefn{\theta}^{(0)}(\tv) = a_{0}(\tv)\bfone,
\end{equation}
for some scalar function $a_{0}(\tv)$.  The expansion \eqref{eq:93} then becomes
\begin{equation}
  \label{eq:97}
   \epsilon^{\kappa}A(x_{*})^{T} \baefn{\theta}^{(1)}(\tv ) + \epsilon^{1-\theta}a_{0}'(\tv )\bv(x_{*}) + \bigo(\epsilon)
+ o(\epsilon^{\kappa}) + o(\epsilon^{1-\theta}) = 0,
\end{equation}
where we have used the fact that $\frac{d^{n}}{dx^{n}}A^{T}\bfone = 0$ for all $n\geq0$.
Setting $\kappa = 1$ recovers the outer solution.  The only remaining possibility is to set $\kappa = 1-\theta$, which yields
\begin{equation}
  \label{eq:98}
  \bigo(\epsilon^{1-\theta}):\quad   A(x_{*})^{T} \baefn{\theta}^{(1)}(\tv ) = - a_{0}'(\tv )\bv(x_{*}),
\end{equation}
and since $\bpss(x_{*})^{T}\bv(x_{*}) = 0$, the solution is
\begin{equation}
  \label{eq:99}
  \baefn{\theta}^{(1)}(\tv ) = -a_{0}'(\tv )\bm{\zeta},
\end{equation}
where $\bm{\zeta}$ satisfies \eqref{eq:90}.  Hence, 
\begin{equation}
  \label{eq:100}
  \baefn{\theta}(\tv ) \sim a_{0}(\tv ) \bfone - \epsilon^{1-\theta}a_{0}'(\tv )\bm{\zeta}.
\end{equation}
The function $a_{0}(\tv )$ is determined at higher order; we find
\begin{multline}
\label{eq:102}
 \epsilon^{1-\theta} A(x_{*})^{T} \baefn{\theta}^{(2)}
+ \epsilon^{\theta}a_{0}'\left(z\bv'(x_{*}) - zA'(x_{*})^{T} \bm{\zeta}\right)\\
   - \epsilon^{1-\theta}a_{0}''\left(\diag{\bv(x_{*})} \bm{\zeta} - \bb(x_{*})\right) = 0.
\end{multline}
Setting $\theta=1/2$ yields
\begin{equation}
  \label{eq:103}
\bigo(\epsilon):\quad A(x_{*})^{T} \baefn{\theta}^{(2)}(\tv ) =  a_{0}''(\tv )\diag{\bv(x_{*})} \bm{\zeta} - a_{0}'(\tv )\tv \left(\bv'(x_{*}) - A'(x_{*})^{T} \bm{\zeta}\right) ,
\end{equation}
and the resulting solvability condition is
\begin{equation}
  \label{eq:104}
   a_{0}''(\tv )  - \tv \left(\frac{(\bpss(x_{*})^{T}\bv(x_{*}))' }{\bpss(x_{*})^{T}\left(\diag{\bv(x_{*})} \bm{\zeta} - \bb(x_{*})\right)}\right) a_{0}'(\tv ) = 0.
\end{equation}
Note that $A\bpss=0$ $\Rightarrow A'\bpss=-A\bpss'$. Furthermore, $\bm{\zeta}^{T}A(x_{*})\bpss'(x_{*}) = \bv(x_{*})^{T}\bpss'(x_{*})$.  Hence, $\bpss(x_{*})^{T}\bv'(x_{*}) - \bpss(x_{*})^{T}A'(x_{*})\bm{\zeta} = (\bpss(x_{*})^{T}\bv(x_{*}))' = \bar{v}(x)'$, where $\dot{x} = \bar{v}(x)$ is the deterministic limit \eqref{eq:114}.  One can show that (see Appendix \ref{sec:curvature-prefactor})
\begin{equation}
  \label{eq:105}
\frac{\bar{v}(x)' }{\bpss(x_{*})^{T}(\bb(x_{*})-\diag{\bv(x_{*})} \bm{\zeta} )} = -\Phi''(x_{*}).
\end{equation}
Assuming that $\Phi''(x_{*})<0$,  the solution to \eqref{eq:104} is
\begin{gather}
  \label{eq:106}
  a_{0}'(\tv ) = \hat{c}_{1} e^{\frac{1}{2}\Phi''(x_{*})\tv ^{2}},\\
  a_{0}(\tv ) = \hat{c}_{0} + \hat{c}_{1}\int_{0}^{\tv }e^{\frac{1}{2}\Phi''(x_{*})\tv'^{2}}d\tv',
\end{gather}
where $\hat{c}_{0,1}$ are unknowns constants.  Note that since $x=x_{*}$ is a local maxima of $\Phi(x)$, we assume that $\Phi''(x_{*})<0$ so that $a_{0}(\tv )\to 0$ as $\tv \to\infty$.  The solution \eqref{eq:100} becomes
\begin{equation}
  \label{eq:108}
  \baefn{\theta}(\tv ) \sim  \left(  \hat{c}_{0} + \hat{c}_{1}\int_{0}^{\tv }e^{\frac{1}{2}\Phi''(x_{*})\tv'^{2}}d\tv' \right)\bfone - \epsilon^{1/2}\hat{c}_{1}e^{\frac{1}{2}\Phi''(x_{*})\tv ^{2}} \bm{\zeta},
\end{equation}
which replaces the first two terms in \eqref{eq:91} (i.e., $\hat{c}_{0}\bfone+\hat{c}_{1}(\bm{\zeta} - z\bfone)$).  Notice that the solution is now bounded in the limit $z\to \infty$, which allows us to match it to the outer solution; we require $\lim_{\tv \to\infty}\baefn{\theta}(z)=\bfone$ so that
\begin{equation}
  \label{eq:109}
  \hat{c}_{0} + \hat{c}_{1}\int_{0}^{\infty}e^{\frac{1}{2}\Phi''(x_{*})\tv ^{2}}d\tv  =
 \hat{c}_{0} + \hat{c}_{1}\sqrt{\frac{\pi}{2\abs{\Phi''(x_{*})}}} =1,
\end{equation}
and take
\begin{equation}
  \label{eq:110}
  \hat{c}_{0} = 1-\hat{c}_{1}\sqrt{\frac{\pi}{2\abs{\Phi''(x_{*})}}}.
\end{equation}
As $\tv \to 0$, 
\begin{equation}
  \label{eq:46}
  \baefn{\theta}(\tv ) \sim (\hat{c}_{0}+ \hat{c}_{1}\tv )\bfone - \epsilon^{1/2}\hat{c}_{1}\bm{\zeta}.
\end{equation}
This matches with \eqref{eq:91} if $c_{0}= \hat{c}_{0}$ and $c_{1} = -\epsilon^{1/2}\hat{c}_{1}$.
The remaining unknown constants $c_{j}$, $j = 1,\cdots,M$, are determined using the absorbing boundary condition \eqref{eq:87}, resulting in the linear system of equations,
\begin{equation}
  \label{eq:113}
  \hat{c}_{1}\left( \sqrt{\frac{\pi}{2\abs{\Phi''(x_{*})}}}\bfone + \epsilon^{1/2} \bm{\zeta} \right) 
- \sum_{j=2}^{M}c_{j}\Upsilon_{j} = \bfone.
\end{equation}
\begin{theorem}
\label{pr:aef1}
A uniform asymptotic approximation, valid throughout the boundary layer and transition regions, of the solution to \eqref{eq:86} is given by 
\begin{multline}
\label{eq:111}
  \baefn{}(x) \sim \left[ 1 - \hat{c}_{1}\left(\sqrt{\frac{\pi}{2\abs{\Phi''(x_{*})}}}-\int_{0}^{(x-x_{*})/\epsilon^{1/2}}e^{\frac{1}{2}\Phi''(x_{*})x'^{2}}dx'\right) \right]\bfone \\
- \epsilon^{1/2} \hat{c}_{1}e^{\frac{1}{2}\Phi''(x_{*}) (x-x_{*})^{2}/\epsilon}\bm{\zeta}
 + \sum_{j=2}^{M}c_{j}\Upsilon_{j}e^{-\meig_{j}(x-x_{*})/\epsilon},
\end{multline}
where
\begin{align}
  \label{eq:4}
  \hat{c}_{1} &\sim \sqrt{\frac{2\abs{\Phi''(x_{*})}}{\pi}} - \epsilon^{1/2}\sqrt{\frac{2\abs{\Phi''(x_{*})}}{\pi}}\hat{c}_{1}^{(1)} + O(\epsilon),\\
  c_{j} &\sim -\epsilon^{1/2}c_{j}^{(1)} + O(\epsilon).
\end{align}
The constants $\hat{c}_{1}^{(1)}$ and $c_{j}^{(1)}$, $j = 2,\cdots, M$, satisfy
\begin{equation}
  \label{eq:74}
  \hat{c}_{1}^{(1)}\bfone + \sum_{j=2}^{M}\hat{c}_{j}^{(1)}\Upsilon_{j} = \sqrt{\frac{2\abs{\Phi''(x_{*})}}{\pi}}\bm{\zeta}.
\end{equation}
\end{theorem}

\subsubsection{Discrete process}
\label{sec:discrete-process}
The adjoint eigenfunction for the discrete process satisfies
\begin{equation}
  \label{eq:115}
  \left( \qss [A(x)]^{T} + \diag{\mathbbm{d}^{*}}\right)\baefn{}(x) = 0,
\end{equation}
where
\begin{equation}
  \label{eq:116}
  \mathbbm{d}^{*}(s) = \qd\left[\hbdrate_{-}(x|s)(\bbe^{-\partial x} - 1) + \hbdrate_{+}(x|s)(\bbe^{\partial x} - 1)\right],
\end{equation}
with $\bbe^{\pm \partial x}$ defined by \eqref{eq:55}.  The absorbing boundary condition is $\baefn{}(x_{*}) = 0$.  Once again, the outer solution is $\baefn{\mathrm{out}} = \bfone$.

Motivated by the boundary layer analysis in Section \ref{sec:semi-cont-proc}, we rescale with $x=x_{*}+\epsilon^{\theta}\tv $.  We are interested in two cases: $\theta=1$ and $\theta=1/2$.  In the former case, the scaling simply returns the process to a discrete variable since $x=\frac{n}{\qd} = \hgam \epsilon n$.  Let $\hat{n}=n-n_{*}$ and $\baefn{\rm bl}(\hat{n}) = \baefn{}( x_{*}+\hgam \epsilon\hat{n})$.  Then to leading order
\begin{multline}
  \label{eq:190}
\hgam[A(x_{*})]^{T}\baefn{\rm bl}(\hat{n})  + \diag{\hbdrate_{-}(x_{*} | s)}\left(\baefn{\rm bl}(\hat{n}-1)-\baefn{\rm bl}(\hat{n})\right) \\
+\diag{\hbdrate_{+}(x_{*} | s)}\left(\baefn{\rm bl}(\hat{n}+1)-\baefn{\rm bl}(\hat{n})\right) = 0.
\end{multline}
Solutions have the form $\baefn{\rm bl}(\hat{n}) = \Gamma_{j}\mu_{j}^{\hat{n}}$.
Substituting this into \eqref{eq:190} yields 
\begin{equation}
  \label{eq:156}
  \left[\hgam \mu_{j}[A(x_{*})]^{T} + \mu_{j}(\mu_{j}-1)\diag{\hbdrate_{+}(x_{*} | s)} - (\mu_{j}-1)\diag{\hbdrate_{-}(x_{*} | s)}\right]\Gamma_{j}=0.
\end{equation}
As before (see \eqref{eq:91}), one of the linearly independent solutions is
\begin{equation}
  \label{eq:191}
  \baefn{\rm bl}(\hat{n}) = \bm{\zeta} - \hgam \hat{n} \bfone,
\end{equation}
where $\bm{\zeta}$ is given by \eqref{eq:90}.  
On the other hand, if $\theta=1/2$ we recover \eqref{eq:93}, which means that we can replace \eqref{eq:191} by \eqref{eq:108}.
We assume that $\abs{\mu_{j}}<1$, $j = 2,\cdots,M$.  
The boundary condition $\baefn{}(x_{*}) = 0$,  results in a linear system having the same form as \eqref{eq:113}.
\begin{theorem}
\label{pr:aef2}
A uniform asymptotic approximation, valid throughout the boundary layer and transition regions, of the solution to \eqref{eq:115} is given by 
\begin{multline}
  \label{eq:195}
    \baefn{}(x) \sim \left[ 1 - \hat{c}_{1}\left(\sqrt{\frac{\pi}{2\abs{\Phi''(x_{*})}}}-\int_{0}^{(x-x_{*})/\epsilon^{1/2}}e^{\frac{1}{2}\Phi''(x_{*})x'^{2}}dx'\right) \right]\bfone \\
- \epsilon^{1/2} \hat{c}_{1}e^{\frac{1}{2}\Phi''(x_{*}) \frac{(x-x_{*})^{2}}{\epsilon}}\bm{\zeta}
 + \sum_{j=2}^{M}c_{j}\Gamma_{j}\mu_{j}^{\frac{x-x_{*}}{\hgam\epsilon}},
\end{multline}
where $\Gamma_{j}$ and $\mu_{j}$ satisfy \eqref{eq:156}.    The constants $\hat{c}_{1}$ and $c_{j}$, $j=2, \cdots, M$, are given by \eqref{eq:4} and \eqref{eq:74}, after substituting $\Gamma_{j}$ for $\Upsilon_{j}$.
\end{theorem}

\subsection{Principal eigenvalue}
\newcommand{\delx}{\Delta x}
Now that we have approximations for the right and left eigenfunction, we can construct the approximation of the principal eigenvalue using the spectral projection method (see \eqref{eq:38}) outlined in the introduction of this section.
\begin{theorem}
\label{pr:lambda}
Let $\eigabs{-}$ and $\eigabs{+}$ be defined for the domain $x<x_{*}$ and $x>x_{*}$, respectively.
Given the asymptotic approximation of the eigenfunction in Theorem \ref{pr:ef} and the adjoint eigenfunction in Theorem \ref{pr:aef1} or \ref{pr:aef2}, an asymptotic approximation of the principal eigenvalue is
\begin{equation}
  \label{eq:155}
  \eigabs{\pm} \sim \left(\frac{B}{\pi}\sqrt{\abs{\Phi''(x_{*})} \Phi''(x_{\pm})}\right)\frac{k(x_{*})}{k(x_{\pm})}\explr{-\frac{1}{\epsilon}(\Phi(x_{*})-\Phi(x_{\pm}))},
\end{equation}
where
\begin{equation}
  \label{eq:77}
  B = \sum_{s}\rho(s | x_{*})\left(b(s, x_{*}) - v(s, x_{*})\zeta(s)\right)
\end{equation}
and
\begin{equation}
\label{eq:64}
  \frac{k(x_{*})}{k(x_{\pm})} = \explr{-(\Psi(x_{*}) - \Psi(x_{\pm}))},
\end{equation}
with $\Phi$, $\Psi$, $\rho$, and $\zeta$ defined by \eqref{eq:58}, \eqref{eq:242}, \eqref{eq:19}, and \eqref{eq:90}, respectively.
\end{theorem}

For both processes, the normalization constant (the denominator in \eqref{eq:38}) is approximated using Laplace's method with
\begin{equation}
  \label{eq:163}
  \ave{\befn{},\baefn{}}\sim \left(\frac{\Phi''(x_{\pm})}{2\pi \epsilon}\right)^{-1/2}.
\end{equation}
The boundary contribution (the numerator in \eqref{eq:38}) for each process is computed as follows.

First, for the semi-continuous process, substituting the eigenfunctions \eqref{eq:63} and \eqref{eq:111} into \eqref{eq:204} yields
\begin{equation}
  \label{eq:73}
  J(\befn{},\baefn{}) \sim  B\sqrt{\frac{2\epsilon\abs{\Phi''(x_{*})}}{\pi}}k(x_{*}) e^{-\frac{1}{\epsilon}\Phi(x_{*})},
\end{equation}
where
\begin{equation}
  \label{eq:56}
  B \equiv \bpss(x_{*})^{T}\bb(x_{*}) - \sqrt{\frac{\pi}{2\abs{\Phi''(x_{*})}}}\sum_{j=2}^{M}\hat{c}_{j}^{(1)}\gamma_{j}\bpss(x_{*})^{T}\diag{\bb(x_{*})}\Upsilon_{j}.
\end{equation}
From \eqref{eq:89} we have that $\meig_{j}\bpss(x_{*})^{T}\diag{\bb(x_{*})}\Upsilon_{j} = \bpss(x_{*})^{T}\diag{\bv(x_{*})}\Upsilon_{j}$.
Then, using \eqref{eq:74} and $\bpss(x_{*})^{T}\bv(x_{*}) = 0$, it follows that \eqref{eq:56} can be rewritten as \eqref{eq:77}.

The discrete version of \eqref{eq:40} can be obtained using a summation by parts argument.  The resulting boundary contribution is
\begin{equation}
  \label{eq:76}
      J(\befn{},\baefn{})  = \frac{1}{2}\befn{}(x_{*})^{T}\left(\diag{\hbdrate_{+}(x_{*} | s)}\baefn{}(x_{*} + \frac{1}{\qd}) - \diag{\hbdrate_{-}(x_{*} | s)} \baefn{}(x_{*} - \frac{1}{\qd})\right).
\end{equation}
The first two terms in $\baefn{}(x_{*}\pm 1/\qd)$ (see \eqref{eq:195}) can be expanded in $1/\qd \ll 1$ (the third term is the boundary layer solution).
Substituting \eqref{eq:195} and \eqref{eq:63} into \eqref{eq:76} (using \eqref{eq:117}) shows that $J$ takes the form \eqref{eq:73} with
\begin{multline}
  \label{eq:75}
  B = \bpss(x_{*})^{T}\left(\bb(x_{*}) - \frac{1}{2}\diag{\bv(x_{*})}\bm{\zeta}\right)  \\
+ \frac{1}{2\hat{c}_{1}^{(0)}} \sum_{j=2}^{M}\frac{\hat{c}_{j}}{\mu_{j}}\bpss(x_{*})^{T}(\mu_{j}^{2}\diag{\hbdrate_{+}(x_{*} | s)} - \diag{\hbdrate_{-}(x_{*} | s)})\Gamma_{j}.
\end{multline}
From \eqref{eq:156} we have that $\mu_{j}\bpss(x_{*})^{T}\diag{\hbdrate_{+}(x_{*} | s)}\Gamma_{j} = \bpss(x_{*})^{T}\diag{\hbdrate_{-}(x_{*} | s)}\Gamma_{j}$, so that
\begin{equation*}
  \frac{1}{\mu_{j}}\bpss(x_{*})^{T}\left(\mu_{j}^{2}\diag{\hbdrate_{+}(x_{*} | s)} - \diag{\hbdrate_{-}(x_{*} | s)}\right)\Gamma_{j} = -\bpss(x_{*})^{T}\diag{\bv(x_{*})}\Gamma_{j}.
\end{equation*}
Thus, we can rewrite \eqref{eq:75} as \eqref{eq:77}.

\section{Example: stochastic model of gene expression}
\label{sec:model}
Consider the following as an example of a discrete Markov process with an internal state.
A population of proteins is modeled as a birth/death process, where the protein production rate depends on the internal state.
The hypothetical gene responsible for producing the protein is said to be activated if an activator molecule is bound to the gene's promotor.
When the gene is activated, protein is produced at a higher rate than when it is unactivated.
For simplicity we refer to ``activated'' and ``unactivated as ``on'' and ``off,'' respectively.
All parameters are presented in nondimensional form (see \citep{kepler01a} for the original dimensional version).
The following state diagram, where $N_{n}$ is the state where $n$ proteins are present in the system, represents the external state transitions:
\begin{equation}
  \label{eq:3}
    N_{0} \chem{\synthrate (S(t))}{\delta} 
    N_{1}\chem{\synthrate (S(t))}{2\delta} 
    N_{2}\cdots\chem{\synthrate (S(t))}{n\delta} 
    N_{n}\chem{\synthrate (S(t))}{(n+1)\delta}\cdots,
\end{equation}
where we set $\delta=1$.
The two state stochastic process, $S(t)$, represents the on/off state of the gene; $S(t)=1$ when the gene is on and $S(t)=0$ when it is off. 
The production rate is a function of the gene state, with $\synthrate (0)=\anot \qd$, $\synthrate (1) = \qd$.
The nondimensional parameter $\anot $ controls how much spontaneous protein production occurs when the gene is off, and we assume $0<\anot <1$ so that protein production is higher when the gene is on. 
To get nonlinear phenomena, the internal state transitions must depend on the external state.
Assume that the activator molecule is a dimer of the protein product so that the protein activates its own gene.
A simple model of the gene is given by
\begin{equation}
  \label{eq:5}
  (\mbox{off})\chem{\qss N(t)^{2}/\qd^{2}}{\qss \beta}(\mbox{on}),
\end{equation}
where $N(t)$ is the number of protein copies.
The transition rate matrix and quasi-steady-state distribution are given by
\begin{equation}
 \label{eq:47}
  A(x) \equiv
\left[  \begin{smallmatrix}
    -x^{2}&\beta\\
   x^{2} & -\beta
  \end{smallmatrix}\right],\quad
  \bpss(x) = \bvec{\frac{\beta}{\beta+x^{2}}}{\frac{x^{2}}{\beta + x^{2}}},
\end{equation}
respectively, where $x = n/\qd$. 
The transitions between the two gene states are assumed to be fast by specifying that $1/\qss\ll 1$.
Writing $\rmp_{s}(n, t) = \rmp(s, n, t)$, the CK equation \eqref{eq:53} is
\begin{subequations}
  \label{eq:6}
  \begin{align}
    \pd{}{t}\mathrm{p}_{0}(n,t) &=   \left[(\mathbb{E}^{+}-1)n+\qd\anot  (\mathbb{E}^{-}-1)\right]\mathrm{p}_{0} + \qss\left( - \frac{n^{2}}{\qd^{2}}\mathrm{p}_{0} + \beta \mathrm{p}_{1}\right)  \\
    \pd{}{t}\mathrm{p}_{1}(n,t) &= \left[(\mathbb{E}^{+}-1)n + \qd (\mathbb{E}^{-}-1)\right]\mathrm{p}_{1}
           +\qss\left(\frac{n^{2}}{\qd^{2}}\mathrm{p}_{0} - \beta \mathrm{p}_{1}\right) ,
  \end{align}
\end{subequations}
where the jump operators $\mathbb{E}^{\pm}$ are defined by $\mathbb{E}^{\pm}f(n) = f(n \pm 1)$.

A semi-continuous process \eqref{eq:66} is given by applying a diffusion approximation to \eqref{eq:6}.  The mean number of proteins when the gene is on is $\qd$.  When $\qd\gg 1$, we can rescale to a continuous variable $X(t) = N(t)/\qd$.
It is straight forward to show that the drift in each state is $v(s, x)$, where
\begin{equation}
  \label{eq:10}
  v(0, x) =\anot - x,\quad v(1, x) = 1 - x,
\end{equation}
and the diffusivity is $\epsilon b(s, x)$, where
\begin{equation}
  \label{eq:11}
  b(0, x) = \frac{\hgam}{2}(\anot + x),\quad b(1, x) = \frac{\hgam}{2}(1 + x).
\end{equation}
Recall that in Section \ref{sec:two-classes-markov} we defined the small parameter $\epsilon = 1/\qss = 1/(\hgam\qd)$.
The corresponding CK equation \eqref{eq:52} is
\begin{subequations}
    \label{eq:12}
  \begin{align}
    \pd{}{t}\mathrm{p}_{0}(x,t) &=  -\pd{}{x}[(\anot -x)\mathrm{p}_{0}] + \frac{1}{2\qd}\pdd{}{x}[(\anot +x)\mathrm{p}_{0}] - \qss\left(  x^{2}\mathrm{p}_{0} - \beta \mathrm{p}_{1}\right) \\
    \pd{}{t}\mathrm{p}_{1}(x,t) &= -\pd{}{x}[(1-x)\mathrm{p}_{1}] + \frac{1}{2\qd}\pdd{}{x}[(1+x)\mathrm{p}_{1}] + \qss\left(x^{2}\mathrm{p}_{0} - \beta \mathrm{p}_{1}\right).
  \end{align}
\end{subequations}

\subsection{Quasi-stationary analysis of the example problem}
We now apply the QSA from Section \ref{sec:long-time-asymptotic} to the example problem.
From \eqref{eq:58}, the equation for $\Phi'$ can be expressed as $\mathcal{H}(x, \Phi'(x)) = 0$.  In particular, for the discrete process we have
\begin{equation}
  \label{eq:79}
\begin{split}
  \mathcal{H}_{\mathrm{disc}}(x,\pvar) = \frac{x^{2}}{\hgam^{2}}(e^{-\hgam\pvar}-1)^{2}    + \frac{\anot}{\hgam^{2}} (e^{\hgam\pvar}-1)^{2}
+ \frac{x(\anot+1)}{\hgam^{2}}(e^{-\hgam\pvar}-1)(e^{\hgam\pvar}-1) &\\
  -  \frac{x}{\hgam}(\beta + x^{2})(e^{-\hgam\pvar}-1) 
-  \frac{1}{\hgam}(\beta \anot + x^{2})(e^{\hgam\pvar}-1), 
\end{split}
\end{equation}
and for the semi-continuous process
\begin{equation}
  \label{eq:50}
  \begin{split}
    \mathcal{H}_{\mathrm{sc}}(x,\pvar) = b(0, x)b(1, x)\pvar^{4}  + (b(0, x)v(1, x)+b(1, x)v(0, x))\pvar^{3}&\\
 +  (v(0, x)v(1, x) -(\beta b(0, x) + x^{2}(x)b(1, x)))\pvar^{2}&\\
- (\beta v(0, x) + x^{2}v(1, x) )\pvar.&
\end{split}
\end{equation}
For the discrete problem, $\mathcal{H}_{\mathrm{disc}}(x, \pvar)$ can be transformed to a cubic polynomial in $q = e^{\hgam\pvar}$.   Then, the solutions are given by the positive real roots of
\begin{multline}
  \frac{\hgam^{2}q^{2}}{(q-1)}\mathcal{H}_{\rm disc}(x, \frac{\ln(q)}{\hgam})  = \anot q^{3} - (x+\hgam x^{2} + \anot(1+x+\hgam\beta))q^{2} \\
           +x(1+x+\hgam(\beta+x^{2}) + \anot ) q - x^{2}.
\end{multline}
All of the roots are real, but only one satisfies $q = 1$ at the deterministic fixed points.
Likewise, there is a single suitable root of $\mathcal{H}_{\mathrm{sc}}=0$.
Once $\Phi'$ is calculated, the potential function $\Phi(x)$ is computed numerically by quadrature.\footnote{In practice, we find that the best way of numerically integrating $\Phi'(x)$ and $\Psi'(x)$ is to use Chebychev approximation methods (we use the GNU Scientific Library).}

The pre exponential factor is calculated using $\Psi'$ from \eqref{eq:242} with
\begin{equation}
  \label{eq:9}
  \isd(x) \equiv \bvec{\frac{-h(1, x, \Phi'(x))}{h(0, x, \Phi'(x))-h(1, x, \Phi'(x))}}{\frac{h(0, x, \Phi'(x))}{h(0, x, \Phi'(x)) - h(1, x, \Phi'(x))}},
\quad  \bm{l}(x) =  \bfone - \bvec{\frac{h(1, x, \Phi'(x))}{\beta + x^{2}}}{\frac{h(0, x, \Phi'(x))}{\beta + x^{2}}},
\end{equation}
where
\begin{gather}
  h_{\mathrm{disc}}(s, x, \pvar) \equiv \frac{1}{\hgam}(e^{\hgam\pvar }-1)\left(v(s, 0) - xe^{-\hgam\pvar }\right), \\
  h_{\mathrm{sc}}(s, x, \pvar ) \equiv \pvar v(s, x) + \pvar ^{2}b(s, x).
\end{gather}
The asymptotic approximation of the quasi-stationary density is then given by Theorem \ref{pr:ef}, and the eigenvalue approximation is given by Theorem \ref{pr:lambda}, with
\begin{equation}
  \bm{\zeta} = \frac{1-x_{*}}{2\beta}  \bvec{ \beta+1}{\beta-1},\quad   B = \hgam x_{*} + \frac{(x_{*}-\anot)(1-x_{*})}{\beta+x_{*}^{2}}.
\end{equation}

\subsection{QSS diffusion approximation}
\label{sec:diff-limit-cont}
If jumps in $S(t)$ are much more frequent than jumps in $X(t)$, then $S(t)$ is approximately stationary (i.e., distributed according to the quasi-steady-state distribution conditioned on a fixed value of $X(t)$).  The combined process $(S(t), X(t))$, can be approximated by averaging out $S(t)$ to obtain a Markov process that approximates $X(t)$.   In other words, although $X(t)$ is not Markovian due to its dependence on $S(t)$, it may be approximately Markovian.

A projection method results in a scalar Fokker--Planck equation for the marginal external-state probability density function
\begin{equation}
  \label{eq:24}
  u(x,t) = \mathrm{p}_{0}(x,t) + \mathrm{p}_{1}(x,t).
\end{equation}
For a general discussion of the QSS projection method see \citep{gardiner83a,thomas12a}.  For brevity we only quote the result here (see \citep{kepler01a} for further details).  The result is
\begin{equation}
  \label{eq:25}
  \pd{u}{t} = -\pd{}{x}(\mv{x}u)+\frac{1}{\qss}\pdd{}{x}\left(\ave{\bb(x)}u\right)+\frac{1}{\qss} \pd{}{x}\left(\difu(x)\pd{u}{x}\right),
\end{equation}
where
\begin{gather}
  \label{eq:29}
   \mv{x} = \frac{\beta (\anot-x)}{\beta+x^{2}} + \frac{x^{2}(1-x)}{\beta+x^{2}},\\
   \label{eq:30}
  \ave{\bb(x)} = \frac{\qss}{2\qd}\left(\frac{\beta(\anot+ x)}{\beta+x^{2}} + \frac{x^{2}(1+x)}{\beta+x^{2}}\right),\\
  \label{eq:31}
  \begin{split}
    \difu(x) &= \frac{\beta((\anot-x)-\mv{x})(\anot-x)}{(\beta+x^{2})^{2}} 
                  + \frac{x^{2}((1-x)-\mv{x})(1-x)}{(\beta+x^{2})^{2}}.
  \end{split}
\end{gather}
Protein fluctuations are captured by $\ave{\bb(x)}$ and gene fluctuations by $\difu(x)$.
Define the combined diffusivity
\begin{equation}
  \label{eq:85}
  B(x) = \frac{1}{\qss}\left(\ave{\bb(x)} +  D(x)\right).
\end{equation}
To make comparisons to other approximations, we define
\begin{equation}
  \label{eq:175}
  \Phi'(x) = \frac{-\mv{x}}{\ave{\bb(x)} + \difu(x)},\quad
  \Psi'(x) = \frac{\frac{d}{dx}\ave{\bb(x)}}{\ave{\bb(x)} + \difu(x)}
\end{equation}
The mean escape time to reach $x_{*}$ having starting at one of the stable fixed points, $x_{\pm}$, can be approximated \citep{gardiner83a} by $T_{\pm} \sim 1/\eigabs{\pm}$, where 
\begin{equation}
  \label{eq:153}
  \eigabs{\pm} =\left(\frac{B(x_{*})}{\pi}\sqrt{\abs{\Phi''(x_{*})} \Phi''(x_{\pm})}\right)e^{-\Psi(x_{*})}\explr{-\frac{1}{\epsilon}\Phi(x_{*})}.
\end{equation}

\subsection{Limiting processes}
\label{sec:limiting-processes}
The full model is given by a discrete process that is valid for any value of $\qss >0$ and $\qd>0$.
The semi-continuous process is an approximation of the discrete process if $\qd\gg 1$, and it is valid for any value of $\qss >0$.
If we assume that $\qss$ is also a large parameter then further reduction is possible using a QSS diffusion approximation, call it the QSS process, presented in Section \ref{sec:diff-limit-cont}.  All three versions contain terms that depend on $\qss$ and $\qd$, and if these parameters are assumed to be large, all three should account for contributions of noise in the internal and external state.

Of course, further model reduction is possible by removing one source of noise: either $\qss\to\infty$ or $\qd\to\infty$.  The former is known in the literature as the adiabatic limit (see \citep{kepler01a,walczak05a}), and the later we call the quasi-deterministic (QD) limit.  If both limits are taken, a deterministic dynamical system is obtained.  Note that all three versions of the example problem---the discrete, semi-continuous, and QSS processes---converge to the same deterministic limit \eqref{eq:14}.
However, as we show is in this section, the three approximations do not necessarily converge in the adiabatic limit ($\qd\to\infty$) or the QD limit ($\qss\to\infty$).
In the rest of this section we explore each limit in turn.

\subsubsection{Deterministic limit $\qss\to\infty$ and $\qd\to\infty$}
\label{sec:bist-dynam-determ}
If we take the limit $\qd\to \infty$ and $\qss\to \infty$, the resulting deterministic system \eqref{eq:114} is
\begin{equation}
  \label{eq:14} 
  \dot{x} = \bar{v}(x) = \frac{\beta(\anot - x) + x^{2}(1-x)}{\beta + x^{2}}.
\end{equation}
Assuming that $\anot \ll 1$, the system is described as follows.  For $\beta_{-}<\beta<\beta_{+}$, where $\beta_{-}\sim  4\anot+\bigo(\anot^{2})$ and $\beta_{+}\sim\frac{1}{4}+\frac{\anot}{2}+\bigo(\anot^{2})$, the system is bistable, with an unstable fixed point at 
\begin{equation}
  \label{eq:16}
  x_{*}\sim \frac{1}{2}(1-\sqrt{1-4\beta})+\bigo(\anot)
\end{equation}
\begin{figure}[tbp]
  \centering
  \includegraphics[width=8cm]{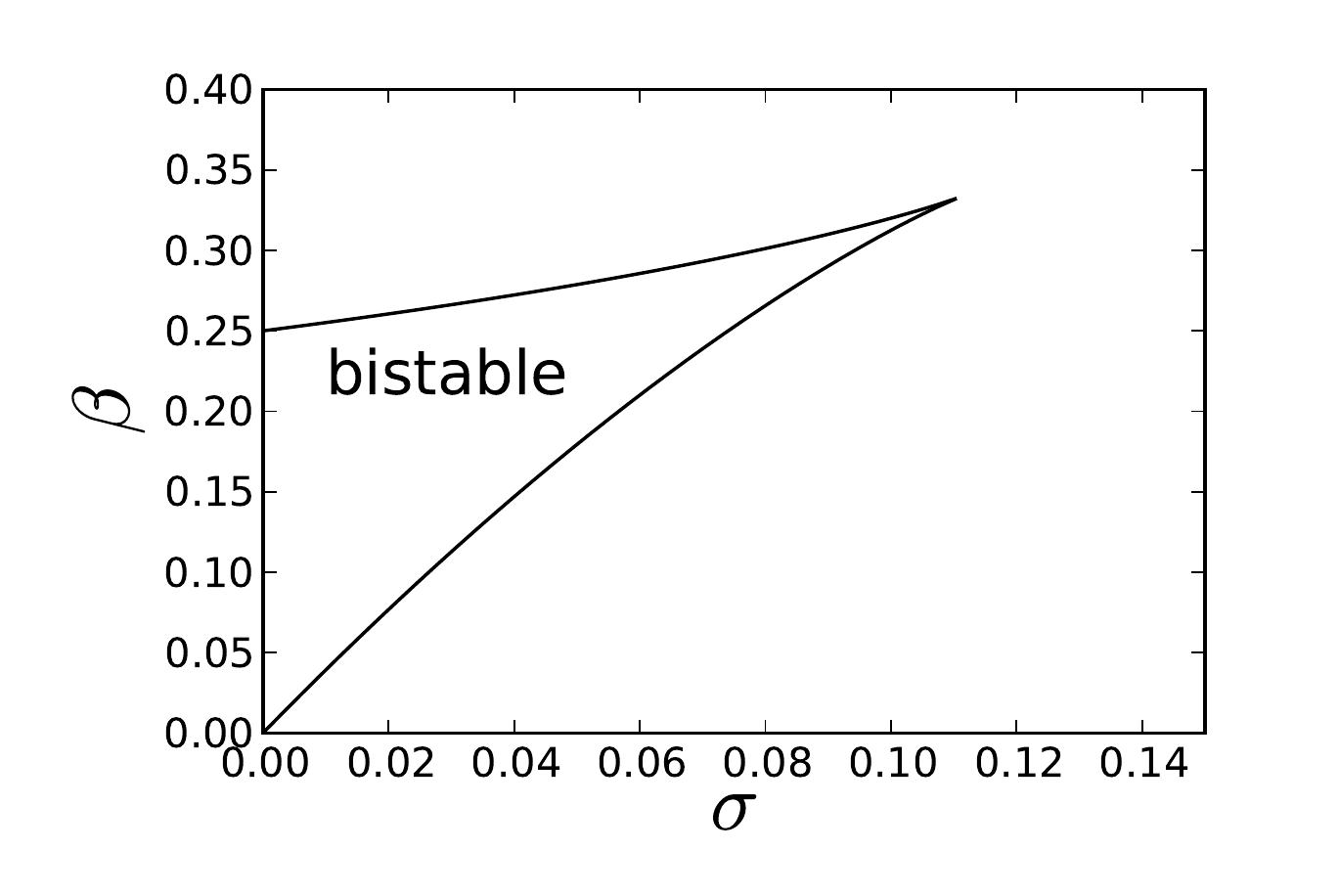}
  \caption{Bifurcation diagram for the deterministic dynamics.}
  \label{fig:bif}
\end{figure}
 and two stable fixed points at
 \begin{equation}
   \label{eq:17}
   x_{-} \sim \anot+\bigo(\anot^{2}), \quad x_{+}\sim \frac{1}{2}(1+\sqrt{1-4\beta})+\bigo(\anot)
 \end{equation}
This is the regime of interest as we wish to characterize the transition times between the two stable fixed points when the system is stochastic with weak fluctuations.

\subsubsection{Quasi-deterministic limit $\qd\to\infty$}
\label{sec:qd}
A velocity jump process can be obtained from the discrete or semi-continuous process by taking the limit $\qd\to\infty$ (both processes converge to the same velocity jump process).  This limit is discussed in \cite{kepler01a} and later a metastable analysis was introduced in \cite{newby12a}.
In this limit, the CK equation converges to
\begin{subequations}
    \label{eq:13}
  \begin{align}
    \pd{}{t}\mathrm{p}_{0}(x,t) &=  -\pd{}{x}[(\anot -x)\mathrm{p}_{0}]  - \qss\left(  x^{2}\mathrm{p}_{0} - \beta \mathrm{p}_{1}\right) \\
    \pd{}{t}\mathrm{p}_{1}(x,t) &= -\pd{}{x}[(1-x)\mathrm{p}_{1}] + \qss\left(x^{2}\mathrm{p}_{0} - \beta \mathrm{p}_{1}\right).
  \end{align}
\end{subequations}
The QSS approximation \eqref{eq:25} does not converge to \eqref{eq:13}; instead the Fokker--Planck equation \eqref{eq:25} becomes
\begin{equation}
  \label{eq:49}
    \pd{u}{t} = -\pd{}{x}(\mv{x}u)+ \frac{1}{\qss} \pd{}{x}\left(\difu(x)\pd{u}{x}\right),
\end{equation}
where $\mv{x}$ and $\difu(x)$ are given by \eqref{eq:29} and \eqref{eq:31}, respectively.

A similar analysis can be carried out on the CK equation \eqref{eq:13} for this process (see \citep{newby11b,keener11a,newby12a} for details), and Theorem \ref{pr:ef} and \ref{pr:lambda} hold.  The result is a fully analytical approximation.  For $\anot < x < 1$ we have
\begin{gather}
  \label{eq:51}
  \Phi(x) =  - \beta\ln(1-x) -\anot^{2}\ln(x-\anot) - \frac{1}{2}(x - \anot)^{2} - 2\anot (x-\anot),\\
  \label{eq:44}
  \isd(x) = \bvec{\frac{1-x}{1-\anot}}{\frac{x-\anot}{1-\anot}},\quad k(x) = \frac{1}{(x - \anot)(1-x)},\quad  B = \frac{(x_{*}-\anot)(1-x_{*})}{\beta+x_{*}^{2}}.
\end{gather}

\subsubsection{Adiabatic limit $\qss\to\infty$}
The CK equation for the semi-continuous process \eqref{eq:12} is asymptotic to \eqref{eq:25} as $\qss \to \infty$.  That is, the semi-continuous process converges to a fully continuous process.
On the other hand , the discrete process converges to a birth/death process as $\qss \to \infty$, which can be derived using a reduction procedure.
The reduction procedure is based on a projection method very similar to the QSS reduction in Section \ref{sec:diff-limit-cont}.  We leave the details to Appendix \ref{sec:a1} and state the result.  The limiting master equation is
\begin{equation}
\label{eq:21}
  \frac{du}{dt} = \qd\mathbb{K}u,
\end{equation}
where $u(n, t) \equiv \sum_{s=0}^{1}\mathrm{p}(s, n, t)$, and
\begin{equation}
  \label{eq:23}
   \mathbb{K} \equiv  (\mathbb{E}^{+}-1)\frac{n}{\qd}   + (\mathbb{E}^{-}-1)f(\frac{n}{\qd}),\quad   f(x) = \frac{\anot \beta + x^{2}}{\beta+x^{2}}.
\end{equation}
Note that at deterministic fixed points, $x_{c}$, we have that $f(x_{c})=x_{c}$.

In the adiabatic limit, we have that $\isd(x) = \bpss(x)$ (see \eqref{eq:20}).  The QSA is well known for the reduced process, and Theorem \ref{pr:ef} holds.  We quote the result here and refer the reader to \citep{schuss10a,doering07a}; that is, $ k(x) = \frac{1}{\sqrt{xf(x)}}$, $B = x_{*}$, and
\begin{gather}
  \label{eq:158}
  \Phi(x) = x\left(\ln(\frac{x}{f(x)})-1\right) 
                      + 2\sqrt{\beta}\left(\tan^{-1}(\frac{x}{\sqrt{\beta}}) - \sqrt{\anot}\tan^{-1}(\frac{x}{\sqrt{\anot\beta}})\right).
\end{gather}
The eigenvalue approximation is
\begin{equation}
  \label{eq:8}
  \eigabs{\pm} \sim \left(\frac{B}{\pi}\sqrt{\abs{\Phi''(x_{*})} \Phi''(x_{\pm})}\right)\frac{k(x_{*})}{k(x_{\pm})}\explr{-\qd(\Phi(x_{*})-\Phi(x_{\pm}))}.
\end{equation}

\section{Results}
\label{sec:results}
\begin{figure}[tbp]
  \centering
  \includegraphics[width=11cm]{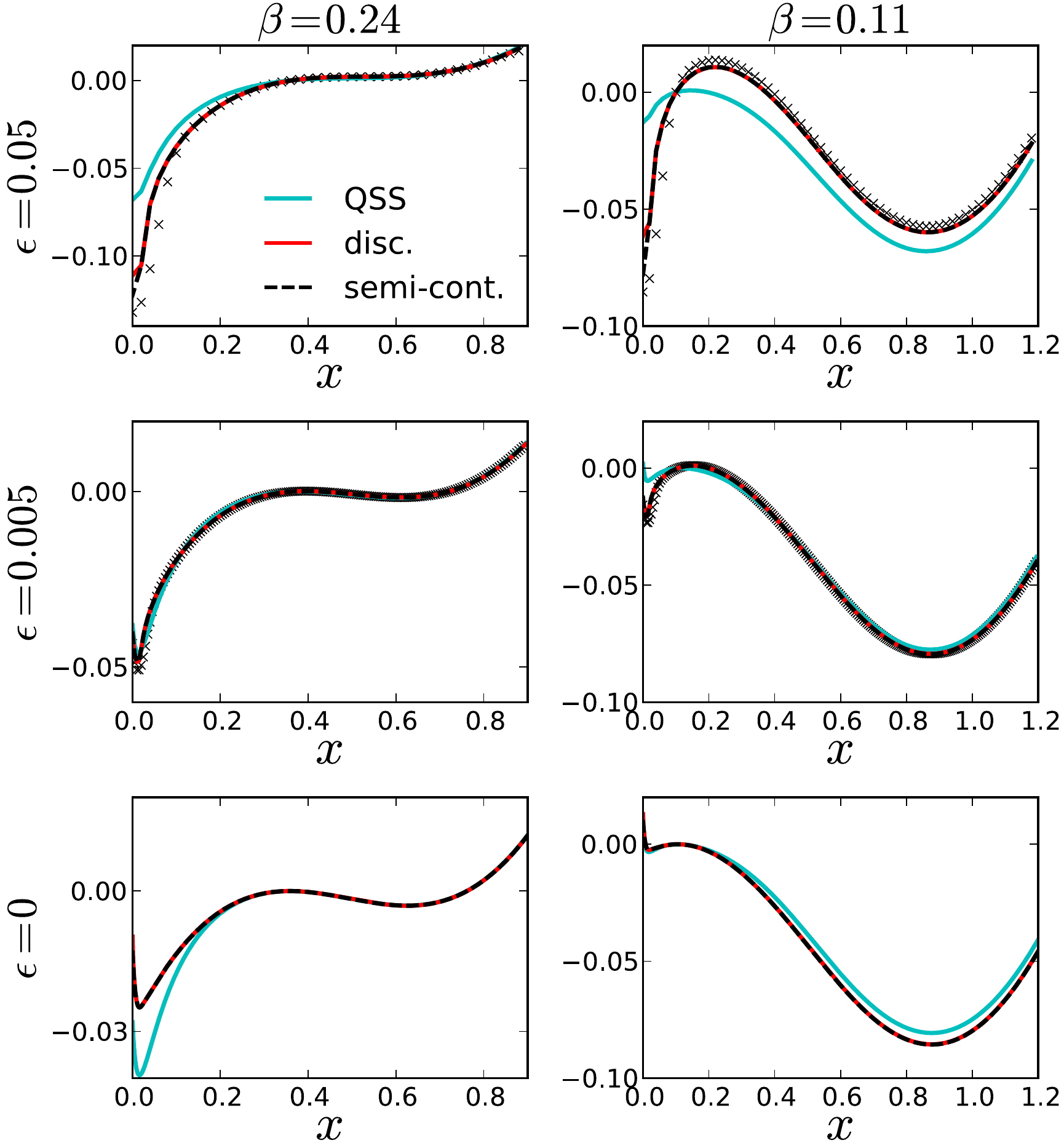}
  \caption{The stability landscape, $-\epsilon\ln(\esd(x))$, for $\hgam=1$ and $\anot = 0.015$.  Each row shows results for a different value of $\epsilon$, and each column shows a different value of the bifurcation parameter $\beta$. The light blue curve is the quasi-steady-state approximation, the green curve is the semi-continuous QSA approximation, and the red curve (which cannot be seen beneath the green curve because both approximations are very close) is the discrete QSA approximation.
For nonzero $\epsilon$ (the top four panes), the first and second order approximation $\Phi(x) + \epsilon\Psi(x)$ is compared to the value of $-\epsilon\ln(p_{\mathrm{s}})$ obtained by a numerical SVD decomposition, shown as ``x'' symbols.  For $\epsilon = 0$ (the bottom two panes), the leading order approximation $\Phi(x)$ is shown.  Note that the SVD solution can not be computed in the $\epsilon \to 0$ limit.}
 \label{fig:stat}
\end{figure}

In this section, we compare the approximations of the stability landscape, defined as $-\epsilon \ln(\esd(x))$, (see \eqref{eq:41}) and of the mean time of a metastable transition from the minimum of one well to the other.  The shape of the stability landscape can be described as a double-well potential, and in Fig.~\ref{fig:stat}
 it is shown for $\anot=0.015$ and two different values of the bifurcation parameter, $\beta$, located within the region of deterministic bistability (see Section \ref{sec:bist-dynam-determ}).  The stability landscape is shown in two columns of plots, each using different parameter values.  In the left column $\beta = 0.24$, which is near the bifurcation point that eliminates the right stability well, and in the right column $\beta = 0.11$, which is near the bifurcation eliminating the left stability well.  Each row shows a different value of $\epsilon$ with $\hgam \equiv \qss/\qd =  1$ so that both noise sources are present.  Approximations of the stability landscape are given by $\Phi(x) + \epsilon \Psi(x)$, where $\Phi$ and $\Psi$ are defined in Section \ref{sec:wkb}.  Note that the WKB approximation of the discrete process breaks down as $x\to 0$ due to small copy number, requiring a boundary correction (see Appendix \ref{sec:xtozero}).  Each approximation---the QSA discrete and semi-continuous approximations,  and the QSS diffusion approximation---is compared to a numerical approximation obtained by SVD decomposition in the top two rows for which $0<\epsilon\ll 1$.  In the bottom row we take the limit $\epsilon\to0$.  Note that the SVD approximation cannot be computed for this case.  First, we observe that the QSA approximation of the discrete and semi-continuous process are so close that they are indistinguishable for every parameter set. (Indeed, we find this to be the case for all of the results presented in this section).  On the other hand, the QSS diffusion approximation shows significant inaccuracies, particularly in the left stability well.  The most significant aspect of the stability landscape that affects metastable transitions is the height of each well in the $\epsilon\to0$ limit.  Although the QSS diffusion approximation does show some error in right stability well, including the height when $\epsilon=0$, these differences are much less significant than the differences in the left well region.  Even for the left well, the QSS diffusion approximation is not always inaccurate.  Indeed, all of the approximations closely agree when $\epsilon = 0.005$ and $\beta=0.24$ (first column, second row of Fig.~\ref{fig:stat}).  However, for other values of $\epsilon$ (top and bottom row) this is clearly not the case.

\begin{figure}[tbp]
  \centering
  \includegraphics[width=11cm]{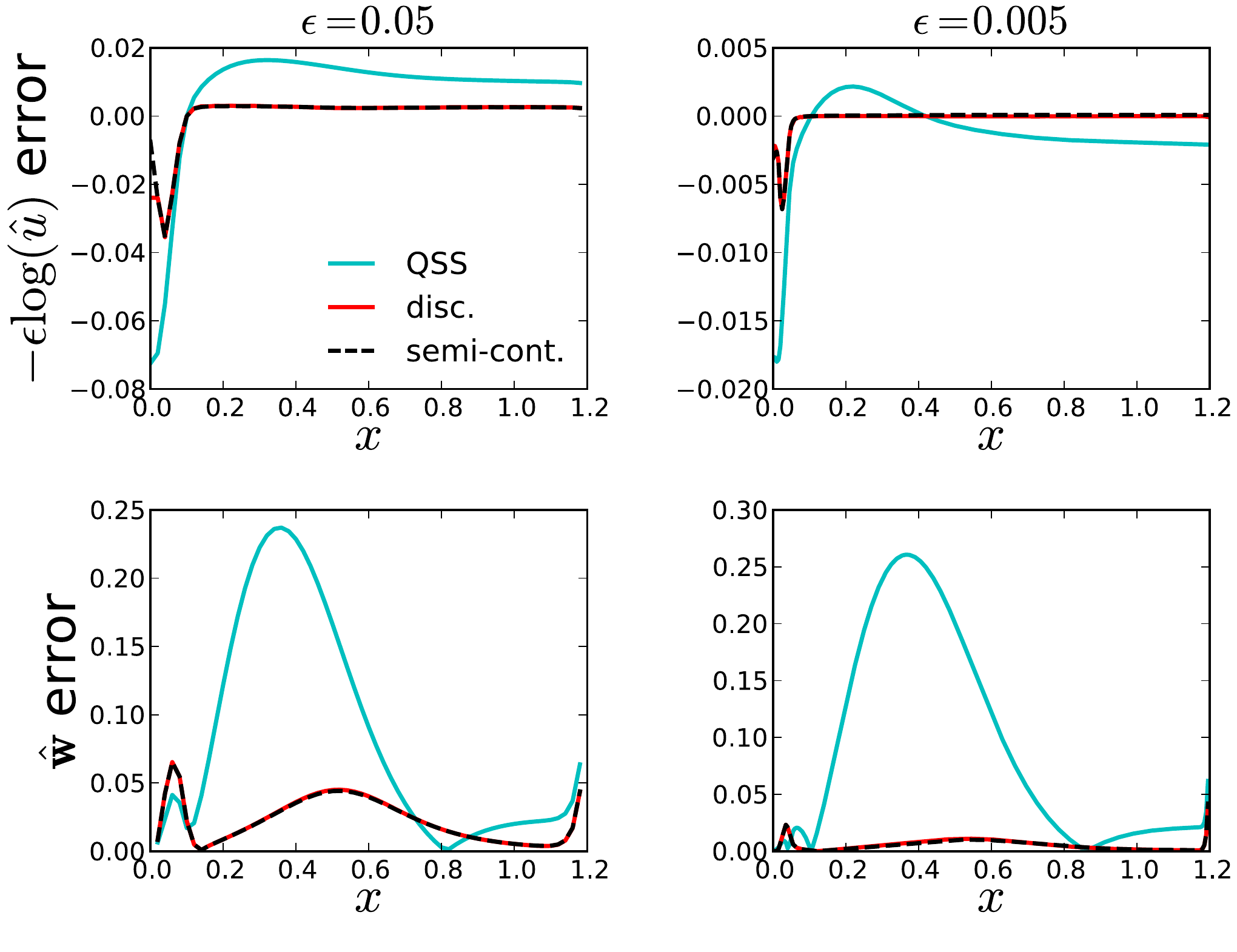}
  \caption{Absolute error.  A comparison of each approximation to the numerical SVD result, for $\beta = 0.11$, $\anot =  0.015$ shown in Fig.~\ref{fig:stat} (right column). The top row shows the error in stability landscape $-\epsilon\ln(\esd(x))$ for $\epsilon = 0.01$ and $\epsilon=0.005$.  The bottom row shows the error in the internal state distribution for the same values of $\epsilon$.  The colors for each curve are the same as in Fig.~\ref{fig:stat}.}
  \label{fig:error}
\end{figure}

To examine the differences in the approximations more closely, we plot the absolute error in the stability landscape and the error in the conditional internal state distribution in Fig.~\ref{fig:error}
 for the parameter values used in the left column of Fig.~\ref{fig:stat} (i.e., $\anot=0.015$, $\hgam=1$, and $\beta=0.11$).  The conditional internal state distribution $\isd(x)$ is \eqref{eq:9} for the discrete and semi-continuous QSA approximations and \eqref{eq:47} for the QSS approximation.  These are again compared to a numerical approximation obtained using an SVD decomposition, and the error is measured using the 1-norm (i.e., $\sum_{s=0}^{1}\abs{\hat{w}_{\mathrm{svd}}(s, x) - \hat{w}_{\mathrm{approx}}(s, x)}$).  The discrete and semi-continuous QSA approximations of the stability landscape show errors primarily in the left well region, while the QSS approximation also shows some error in the right well.  Interestingly, the conditional internal state distribution error is significant for the QSS approximation, peaking at $25\%$ between $x_{*}$ and $x_{+}$.  We expect this error to be quite small near the deterministic fixed points, where all the approximations agree.  We emphasize as one of the key results of this paper that away from fixed points, the conditional internal state distribution is {\em not} always close to the steady-state distribution as assumed in the QSS approximation method.  This has been shown rigorously for velocity jump processes \citep{newby11b}, for which the QD limit is an example.

\begin{figure}[tbp]
  \centering
  \includegraphics[width=11cm]{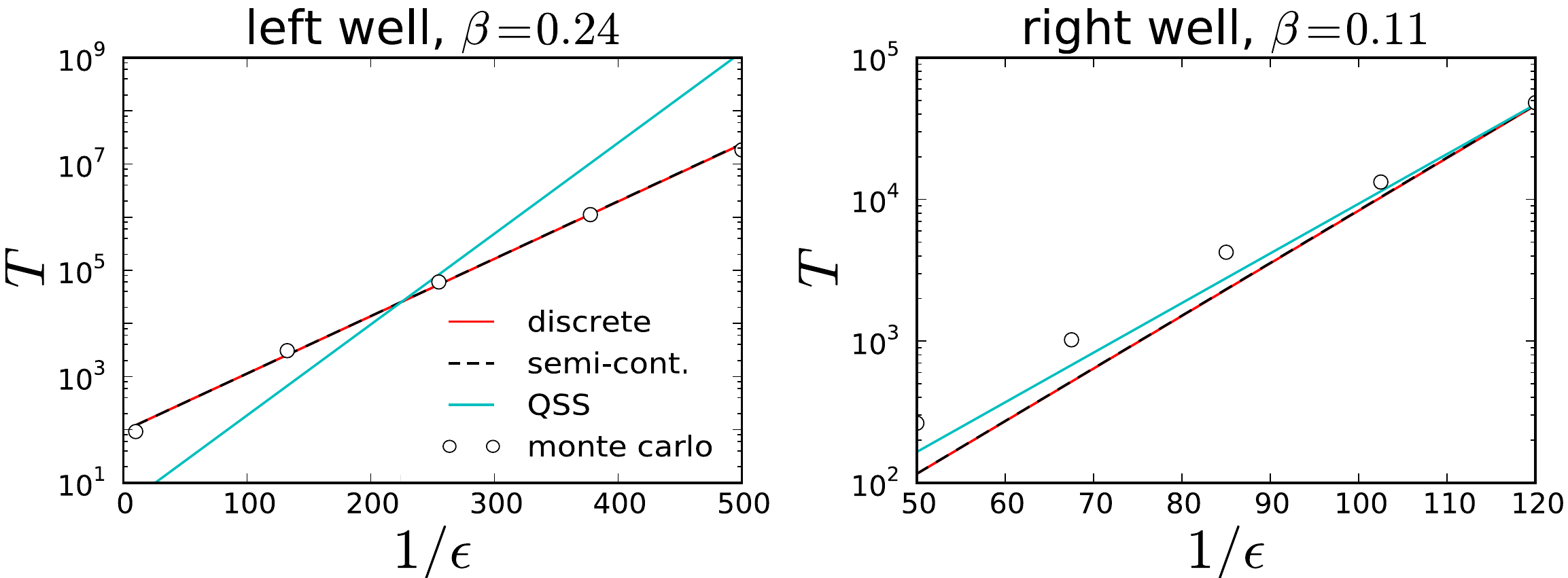}
  \caption{Mean exit time approximations compared to Monte-Carlo results. (a) Exit from the left well for the same parameters as used in the first column of Fig.~\ref{fig:stat}.  (b) Exit from the right well for the same parameters as used in the second column of Fig.~\ref{fig:stat}. }
  \label{fig:mfpt2}
\end{figure}

The approximation of the mean time for a metastable transition between wells is shown in Fig.~\ref{fig:mfpt2}.
The mean escape time approximations, defined as $T_{\pm}\sim 1/\eigabs{\pm}$ (see \eqref{eq:155}), are compared to exact Monte-Carlo (MC) simulations (using the Gillespie algorithm) for parameter values used in Fig.~\ref{fig:stat}.  The mean escape time is plotted on a log scale as a function of $1/\epsilon$ because $\ln(\eigabs{\pm})$ is a linear function of this quantity, with a slope determined by the height of the potential well in the $\epsilon\to 0$ limit (see Fig.~\ref{fig:stat} bottom row).  Escape from the left well (for $\beta=0.11$, Fig.~\ref{fig:stat} left column) is shown on the left, where the discrete and semi-continuous QSA approximations are in good agreement with MC simulations.  The three approximations converge near $\epsilon=0.005$ consistent with Fig.~\ref{fig:stat} (first column, second row).  

A somewhat unexpected result is obtained for escape from the right well (corresponding to the right column of Fig.~\ref{fig:stat}).  All three approximations are very close, and the QSS approximation is actually more accurate for smaller values of $1/\epsilon$.  The difference in the slope of each approximation is slight (see Fig.~\ref{fig:stat} right column, bottom row) and the error in the QSS approximation should grow as $1/\epsilon\to \infty$.  We cannot offer a definitive explanation for the accuracy of the QSS approximation for escape from the right well.  One explanation is that the QSS approximation is valid for larger values of $x$, which seems reasonable since it relies on fast transitions between internal states and the rate of transitioning from the inactive to the active internal state is proportional to $x^{2}$.  However, this is inconsistent with the error in the conditional internal state distribution shown in Fig.~\ref{fig:error} (bottom row), which is the key assumption underlying the QSS approximation.

\begin{figure}[tbp]
  \centering
  \includegraphics[width=11cm]{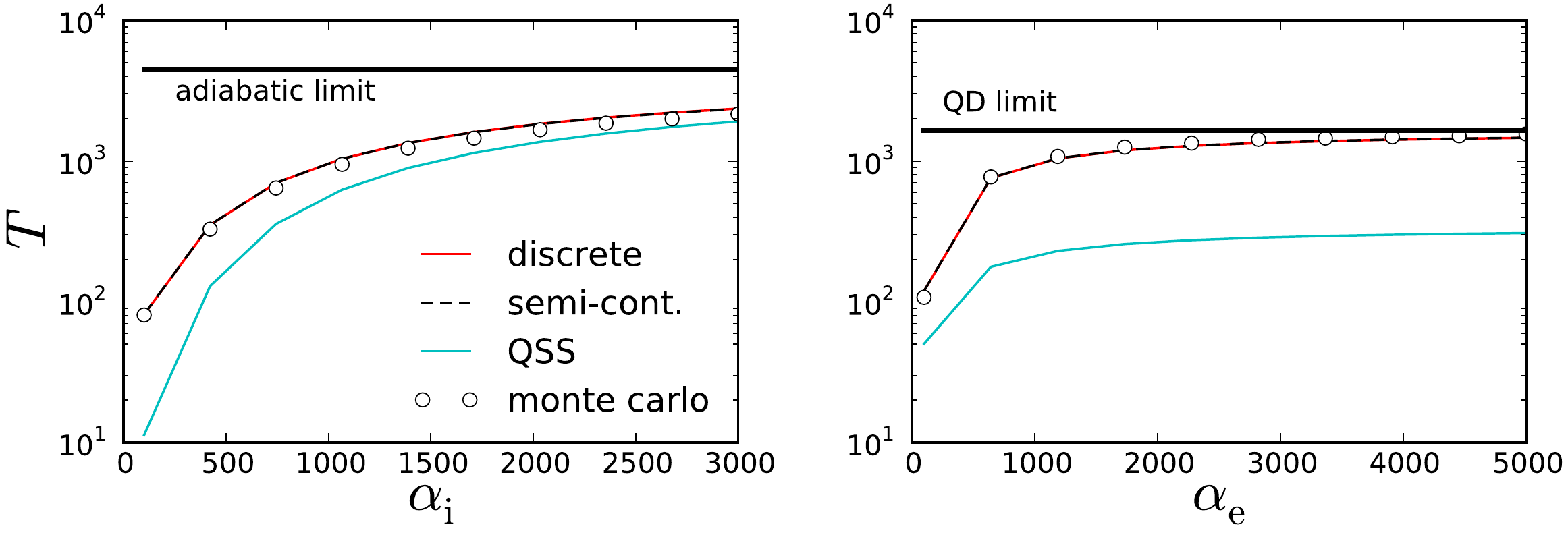}
  \caption{The mean exit time for escape from the left well to the right well, with $\beta = 0.23$, $\anot=0.04$.  Three different approximations (solid curves) are compared to Monte-Carlo simulation results (symbols).  The discrete (red) and semi-continuous (green) QSA approximations are indistinguishable.  Also shown is the QSS approximation (light blue).  (a) The mean exit time as a function of $\qd$ for fixed $\qss = 333$.  (b) The mean exit time as a function of $\qss$ for fixed $\qd=200$.}
  \label{fig:mfpt}
\end{figure}

Finally, we compare the mean escape time in the adiabatic limit $\qss\to\infty$ and in the QD limit $\qd\to\infty$.  In Fig.~\ref{fig:mfpt}, 
the mean time for escape from the left well is shown for $\anot=0.04$ and $\beta=0.23$.  In contrast to previous results, we do not fix $\hgam=\qss/\qd=1$.  Fig.~\ref{fig:mfpt} (right) illustrates that the discrete and semi-continuous approximation converge in the QD limit, and as expected, the QSS approximation error is significant.  In the adiabatic limit (Fig.~\ref{fig:mfpt} left) the discrete and semi-continuous QSA approximations show close agreement for all values of $\qss$, and as expected, all three approximations converge as $\qss\to\infty$.  Even though the discrete and semi-continuous QSA approximations do not converge in the adiabatic limit, the difference is very small.  
This suggests that a diffusion approximation for the external state---recall that we used such a procedure to derive the semi-continuous process from the full discrete process---may be valid in certain situations, which is interesting since diffusion approximations generally break down for metastable behavior due to large deviation errors.  It is possible that the good agreement that we see for the example problem is due to the linear nature of the birth-death process governing transitions in the external state (i.e., that it is due to the simplicity of the example problem).  Since the general QSA procedure presented here does not depend on this assumption, it would be interesting to see how this type of diffusion approximation behaves for a more complicated process.

\appendix
\section{Curvature prefactor}
\label{sec:curvature-prefactor}
The purpose of this section is to show that the part of the eigenvalue estimate that contains information about the curvature of the stability well at the stable and unstable fixed point is unaffected by the QSS diffusion approximation.  This is a reflection of the fact that diffusion approximations, in general, are accurate in a neighborhood of a deterministic fixed point.
The eigenvalue approximation \eqref{eq:155} contains a prefactor term of the form $\sqrt{\abs{\Phi''(x_{*})} \Phi''(x_{\pm})}$.
We would like to show that, when evaluated at a fixed point, $x_{c}$,  the second derivative of $\Phi$ for the discrete, semi-continuous, and QSS processes are all identical.  
We can express the second derivative in terms of $\mathcal{H}$, defined by \eqref{eq:58}, as follows.  

Differentiating $\mathcal{H}(x,\Phi'(x)) = 0$ with respect to $x$ yields
\begin{equation}
  \label{eq:54}
  \frac{d}{dx}\mathcal{H}(x,\Phi'(x)) = \mathcal{H}_{x}(x,\Phi'(x)) + \Phi''(x)\mathcal{H}_{\pvar}(x,\Phi'(x)) = 0,
\end{equation}
and it follows that
\begin{equation}
  \label{eq:57}
  \Phi''(x) = - \frac{\mathcal{H}_{x}(x,\Phi'(x))}{\mathcal{H}_{\pvar}(x,\Phi'(x))}.
\end{equation}
However, we have that
\begin{equation}
  \label{eq:183}
  \mathcal{H}_{\pvar}(x_{c},0) = \mathcal{H}_{x}(x_{c},0) = \mathcal{H}_{xx}(x_{c},0) = 0.
\end{equation}
A formula valid at fixed points can be obtained as follows. Differentiating $\mathcal{H}(x, \Phi'(x)) = 0$ twice with respect to $x$ yields
\begin{equation}
  \label{eq:254}
  \frac{d^{2}}{dx^{2}}\mathcal{H}(x,\Phi'(x)) = \mathcal{H}_{xx}+ \Phi''\mathcal{H}_{x \pvar} 
+  \Phi''(\mathcal{H}_{\pvar x} +  \Phi''\mathcal{H}_{\pvar \pvar}) + \Phi'''\mathcal{H}_{\pvar} = 0,
\end{equation}
and it follows from \eqref{eq:183} that
\begin{equation}
  \label{eq:197}
  \Phi''(x_{c}) = \frac{-2\frac{\partial^{2}}{\partial \pvar \partial x}\mathcal{H}(x_{c},0)}{\pdd{}{\pvar}\mathcal{H}(x_{c},0)}.
\end{equation}
At a fixed point, we have that $p=0$.  Expand $\mathcal{H}(x, p)$ in a Taylors series around $p=0$.  
To second order in $p$, the expansion is consistent with a diffusion approximation, which always corresponds to a Hamiltonian that is quadratic in $p$ with
\begin{equation}
  \label{eq:83}
  \mathcal{H}_{\rm diff}(x, p) = a(x)p + g(x)p^{2},
\end{equation}
where $a(x)$ is the drift and $g(x)$ is the scaled diffusivity.  
For a QSS diffusion approximation of the processes described in Section \ref{sec:two-classes-markov}, one can show that 
\begin{equation}
  \label{eq:84}
  a(x) = \bpss(x)^{T}\bv(x),\quad g(x) = \bpss(x)^{T}\bb(x) - \bpss(x)^{T}\left(\diag{\bv(x)}-a(x)I\right)[A^{\dag}(x)]^{T}\bv(x).
\end{equation}
It follows that at a fixed point $a(x_{c}) = 0$ and $ g(x_{c}) = \bpss(x_{c})^{T}(\bb(x) - \diag{\bv(x_{c})}\bm{\zeta})$.
Substituting \eqref{eq:83} into \eqref{eq:197} yields $\Phi''(x_{c}) = -\frac{a'(x_{c})}{g(x_{c})}$

\section{Adiabatic limit of the discrete process}
\label{sec:a1}
Consider the Master equation for the probability distribution function $\mathrm{p}_{j}(\bk,t) \equiv \mathrm{p}(j,\bk,t | j_{0},\bk_{0},t_{0})$.
In matrix/operator form, the CK equation is
\begin{equation}
\label{eq:123}
  \frac{d\bp}{dt} = L_{1}p + \frac{1}{\epsilon}L_{2}p,
\end{equation}
where $\bp(\bk,t)=(\mathrm{p}_{1}(\bk,t),\: \mathrm{p}_{2}(\bk,t),\cdots,\:\mathrm{p}_{M}(\bk,t))^{T}$; $L_{1}=\diag{\mathcal{D}_{j}}$ is a diagonal matrix of linear operators acting on $\bk$, each of which has a $\bbW$-matrix representation; $L_{2}$ is an $M\times M$  $\bbW$-matrix governing the transitions between internal states, with transition rates that may depend on $\bk$.  Define the projection operator $  \mathcal{P} \equiv \bpss\bfone^{T}$, where $L_{2}\bpss(\bk)=0$, with $\bpss(\bk)>0$ and $\sum_{j=1}^{M}\pss{j}(\bk) = 1$; and $\bfone\equiv (1,\: 1,\cdots,\:1)^{T}$.  We assume the solution has the following form
\begin{equation}
\label{eq:125}
  \bp(\bk,t) = \mathcal{P}\bp(\bk,t) + (I-\mathcal{P})\bp(\bk,t)=u(\bk,t)\bpss(\bk) + \epsilon \bw(\bk,t),
\end{equation}
where
\begin{equation}
\label{eq:126}
  u(\bk,t) \equiv \bfone^{T}\bp(\bk,t),\quad \bfone^{T}\bw(\bk,t) = 0.
\end{equation}
Applying the projection operator to both sides of \eqref{eq:123} yields
\begin{equation}
\label{eq:127}
  \frac{d u}{dt} \bpss = \mathcal{P}L_{1}(u\bpss + \epsilon \bw).
\end{equation}
On the other hand, applying the orthogonal projection yields
\begin{equation}
\label{eq:128}
\epsilon \frac{d\bw}{dt} - \epsilon (I-\mathcal{P})L_{1}\bw = (I-\mathcal{P})L_{1}(u\bpss) + L_{2}\bw.
\end{equation}
After setting $\epsilon=0$ in the above equation we get
\begin{equation}
\label{eq:129}
  \bw(\bk,t) \sim -L_{2}^{-1}(I-\mathcal{P})L_{1}(u(\bk,t)\bpss(\bk)).
\end{equation}
Substituting \eqref{eq:129} into \eqref{eq:127} yields the scalar-valued operator equation for $u(\bk,t)$
\begin{equation}
\label{eq:130}
  \frac{du}{dt} = \bfone^{T}L_{1}(u\bpss) - \epsilon \bfone^{T}L_{1}L_{2}^{-1}(I-\mathcal{P})L_{1}(u\bpss).
\end{equation}
One can rewrite \eqref{eq:130} in matrix form to obtain a linear system of ODEs for the vector $\bu(t)$ with elements $u_{\bk}(t)\equiv u(\bk,t)$
\begin{equation}
\label{eq:131}
  \frac{d\bu}{dt} = W \bu,
\end{equation}
where $  W \equiv \sum_{j=1}^{M}\mathcal{D}_{j}\pss{j}$.
In general, the reduced equation represents a Markov process only at leading order.

\section{WKB/KM expansion}
\label{sec:wkbkm-expansion}
Consider the action of the operator $\bbe^{\partial x}$ on $g(x)e^{-\qd\tilde{\Phi}(x)}$ where $g(x)$ is scalar function and $\tilde{\Phi}(x) = \hgam\Phi(x)$.  We have that
\begin{equation}
  \label{eq:70}
\begin{split}
& \bbe^{\pm \partial x}\left(g(x)e^{-\qd \tilde{\Phi}(x)}\right) \\
&\qquad= \sum_{n=0}^{\infty}\frac{(\pm1)^{n}}{n!\qd^{n}}\frac{d^{n}}{dx^{n}}\left[g(x)e^{-\qd \tilde{\Phi}(x)}\right]\\
&\qquad=  \sum_{n=0}^{\infty}\frac{(\pm1)^{n}}{n!\qd^{n}}\sum_{k=0}^{n}\binom{n}{k}g^{(n-k)}(x)\frac{d^{k}}{dx^{k}}e^{-\qd \tilde{\Phi}(x)}\\
&\qquad=  \sum_{n=0}^{\infty}\frac{(\pm1)^{n}}{n!\qd^{n}}\sum_{k=0}^{n}\binom{n}{k}g^{(n-k)}(x)(1+\mathrm{B}_{k}(-\qd \tilde{\Phi}(x)))e^{-\qd \tilde{\Phi}(x)},
\end{split}
\end{equation}
where $\mathrm{B}_{k}$ is the $k$th complete Bell polynomial,
\begin{equation}
  \label{eq:71}
  \mathrm{B}_{n}(f(x)) \equiv \det
  \begin{bmatrix}
    f' & \binom{n-1}{1}f'' & \binom{n-1}{2} f^{(3)} & \cdots & f^{(n)}\\
    -1 & f'                         & \binom{n-2}{1}f''        & \cdots & f^{(n-1)} \\
     0 & -1                        & f'                                 & \cdots & f^{(n-2)} \\
     \vdots & \vdots                & \ddots                                 & \ddots & \vdots \\
     0 & 0                & \cdots                                 & -1 & f' \\
  \end{bmatrix},
\end{equation}
and $B_{0}=0$.  One can show that
\begin{equation}
  \label{eq:72}
  \mathrm{B}_{k}(-\qd\tilde{\Phi}(x)) = \qd^{k}(- \tilde{\Phi}'(x))^{k} - \qd^{k-1}\frac{k}{2}(k-1)\tilde{\Phi}''(x)(-\tilde{\Phi}'(x))^{k-2} + \bigo(\qd^{k-2}).
\end{equation}
Expanding \eqref{eq:70} in terms of $1/\qd$ yields
\begin{multline}
  \label{eq:43}
\bbe^{\pm\partial x}\left(g(x)e^{-\qd \tilde{\Phi}(x)}\right) \\  
= e^{\mp\tilde{\Phi}'(x)}\left[g(x) - \frac{1}{\qd}\left(g'(x) \mp \frac{1}{2}g(x)\tilde{\Phi}''(x)\right)  + \bigo(\qd^{-2})\right]e^{-\qd \tilde{\Phi}(x)}.
\end{multline}

\section{Evaluating $\lim_{x\to x_{c}}\Psi'(x)$ for $x_{c}=x_{\pm},x_{*}$}
\label{sec:fplim}
Using L'H\^opital's rule, we find that
\begin{equation}
  \label{eq:252}
  \begin{split}
    \Psi'(x_{c}) =  \Bigg[ {H}_{\pvar x x} + \frac{1}{2}\Phi''(x_{c})(3{H}_{\pvar\pvar x} + \Phi''(x_{c}){H}_{\pvar\pvar\pvar}) + \frac{1}{2}\Phi'''(x_{c}){H}_{\pvar \pvar}  \qquad\qquad&\\
                         +\bm{l}'(x_{c})^{T}\mathbf{H}_{\pvar x}(x_{c}, 0) + \frac{1}{2}\Phi''(x_{c})\bm{l}'(x_{c})^{T}\mathbf{H}_{\pvar\pvar}(x_{c}, 0)  \Bigg]\quad& \\
      \Big/ \Big[{\bm{l}'(x_{c})^{T}\mathbf{H}_{\pvar}(x_{c}, 0) + {H}_{\pvar x}+ \Phi''(x_{c}){H}_{\pvar\pvar}}\Big], \quad&
  \end{split}
\end{equation}
where $\mathbf{H}(x, p)$ is defined by \eqref{eq:243} and partial derivatives of $H(x, p) \equiv \bfone^{T}\mathbf{H}(x, p)$ are evaluated at $x=x_{c}$ and $p=0$, as for example,
\begin{equation}
  \label{eq:253}
  H_{xp} \equiv \bfone^{T}\frac{\partial^{2}}{\partial x \partial p}\mathbf{H}(x_{c}, 0).
\end{equation}
We also have that $\Phi''(x_{c})$ is given by \eqref{eq:197}, and
\begin{equation}
  \label{eq:251}
\Phi'''(x_{c}) = -2\frac{\mathcal{H}_{\pvar x x}(x_{c}, 0) + \frac{1}{3}\Phi''(x_{c})\mathcal{H}_{\pvar \pvar \pvar}(x_{c}, 0)}{\mathcal{H}_{\pvar \pvar}(x_{c}, 0)}.
\end{equation}
Note that $H(x, p)\neq \mathcal{H}(x, p)$, where $\mathcal{H}(x, p)$ is the Hamiltonian \eqref{eq:58}.

\section{$x\to0$ limit of the quasi-stationary density}
\label{sec:xtozero}
The WKB approximation \eqref{eq:63} of the discrete process breaks down in the limit $x\to 0$, due to small copy number effects (i.e., fluctuations are on the same order). This fact is not relevant if one is interested only in approximating the mean exit time.  However, we also approximate the effective potential.  Although $\Phi(x)$ is bounded in the limit $x\to0$, $\Psi(x)$ has a logarithmic singularity.
To correct this, we use the discrete master equation \eqref{eq:6} to calculate $\befn{}(0)$, with $\befn{}(0) = - (\qss A(0)-\qd\diag{v(0)})^{-1}\befn{}(\frac{1}{\qd})$.


\end{document}